\documentclass[11pt,a4paper,reqno]{article}

\usepackage[utf8]{inputenc}
\usepackage[T1]{fontenc}
\usepackage{lmodern}
\usepackage{amsmath,amssymb,amsthm,amsfonts,mathtools}
\usepackage{mathrsfs}
\usepackage[top=2.8cm, bottom=2.8cm, left=2.8cm, right=2.8cm]{geometry}
\usepackage{microtype}
\usepackage{enumitem}
\usepackage{xcolor}
\usepackage{booktabs}
\usepackage{tikz}
\usetikzlibrary{arrows.meta, quotes}

\usepackage[colorlinks=true,
            linkcolor=blue!80!black,
            citecolor=green!50!black,
            urlcolor=blue!70!black]{hyperref}
\usepackage[capitalize,nameinlink]{cleveref}

\theoremstyle{plain}
\newtheorem{theorem}{Theorem}[section]
\newtheorem{lemma}[theorem]{Lemma}
\newtheorem{proposition}[theorem]{Proposition}

\theoremstyle{definition}
\newtheorem{definition}[theorem]{Definition}

\theoremstyle{remark}
\newtheorem{remark}[theorem]{Remark}

\numberwithin{equation}{section}

\newcommand{\R}{\mathbb{R}}
\newcommand{\N}{\mathbb{N}}

\DeclarePairedDelimiter{\abs}{\lvert}{\rvert}

\DeclareMathOperator{\sech}{sech}

\newcommand{\dxdy}{\,\mathrm{d}x\mathrm{d}y}

\title{\textbf{Shrinking Limit for the Stationary Nonlinear Schrödinger Equation on Cross-Shaped Infinite Domains}}

\author{
    \textbf{Boris Shakarov}\thanks{
        Institut de Mathématiques de Toulouse; UMR 5219,
        Université de Toulouse, CNRS, UPS,
        F-31062 Toulouse Cedex 9, France.
        Email: \texttt{boris.shakarov@math.univ-toulouse.fr}
    }
}

\date{}

\begin{document}

\maketitle

\begin{abstract}
We study the stationary nonlinear Schrödinger equation on a planar, cross-shaped infinite domain and in the shrinking limit. For any fixed thickness, we establish the existence of solutions through a variational approach. To overcome the loss of compactness inherent to the unbounded geometry and bypass the regularity obstructions caused by the re-entrant corners, we constrain the minimization problem to a cone of symmetric and monotonic functions. We then establish the resolvent-preservation property, ensuring the constrained minimizers are stationary solutions. As the transverse width vanishes, we establish a transverse collapse phenomenon, that is, the two-dimensional solutions flatten across the central junction and become purely one-dimensional along the infinite branches. In this way, we prove that cross solutions converge strongly to a stationary solution on the limiting graph in a suitable sense. Furthermore, we demonstrate that this dimensional reduction recovers both vertex continuity and the Kirchhoff conditions.
\smallskip


\end{abstract}

\vspace{0.5cm}

\section{Introduction}

The study of nonlinear dispersive and elliptic equations on metric graphs has gained considerable attention in recent years \cite{zbMATH07618040}. A key question is how solutions to the stationary nonlinear Schrödinger equation on thin higher-dimensional tubular domains concentrate and reduce to one-dimensional models on metric graphs as the channel width vanishes—the so-called \emph{shrinking limit}. Indeed, metric graphs are typically introduced as idealized models for physical networks, such as wires, mesoscopic systems, or optical waveguides.

In this work, we consider the planar cross-shaped infinite domain defined by the union of two perpendicular infinite strips of thickness $L > 0$:
\begin{equation*}
    D(L) \coloneqq \left( \R \times \left(-\frac{L}{2}, \frac{L}{2}\right) \right) \cup \left( \left(-\frac{L}{2}, \frac{L}{2}\right) \times \R \right) \subset \R^2.
\end{equation*}

\begin{figure}[htbp]
    \centering
    \begin{tikzpicture}[
        scale=1.0,
        >=Stealth,
        boundary/.style={thick, black},
        axis/.style={->, thin, gray!80},
        dim/.style={<->, thin, black!80}
    ]
    
        \def\a{1.1}
        \def\m{3.6}

        \fill[gray!15] 
            (-\m, -\a) -- (-\a, -\a) -- (-\a, -\m) -- 
            ( \a, -\m) -- ( \a, -\a) -- ( \m, -\a) -- 
            ( \m,  \a) -- ( \a,  \a) -- ( \a,  \m) -- 
            (-\a,  \m) -- (-\a,  \a) -- (-\m,  \a) -- cycle;

        \draw[boundary] (-\m,  \a) -- (-\a,  \a) -- (-\a,  \m);
        \draw[boundary] ( \m,  \a) -- ( \a,  \a) -- ( \a,  \m);
        \draw[boundary] (-\m, -\a) -- (-\a, -\a) -- (-\a, -\m);
        \draw[boundary] ( \m, -\a) -- ( \a, -\a) -- ( \a, -\m);

        \draw[axis] (-\m - 0.4, 0) -- (\m + 0.4, 0) node[right, black] {$x$};
        \draw[axis] (0, -\m - 0.4) -- (0, \m + 0.4) node[above, black] {$y$};

        \draw[dim] (2.6, -\a) -- (2.6, \a) node[midway, fill=gray!15, inner sep=2pt] {$L$};
        \draw[dim] (-\a, 2.6) -- (\a, 2.6) node[midway, fill=gray!15, inner sep=2pt] {$L$};

        \node at (-2.3, 0.55) {$D(L)$};
    \end{tikzpicture}
    \caption{The cross-shaped domain $D(L) \subset \R^2$ formed by two perpendicular strips of width $L$.}
    \label{fig:domain_DL}
\end{figure}
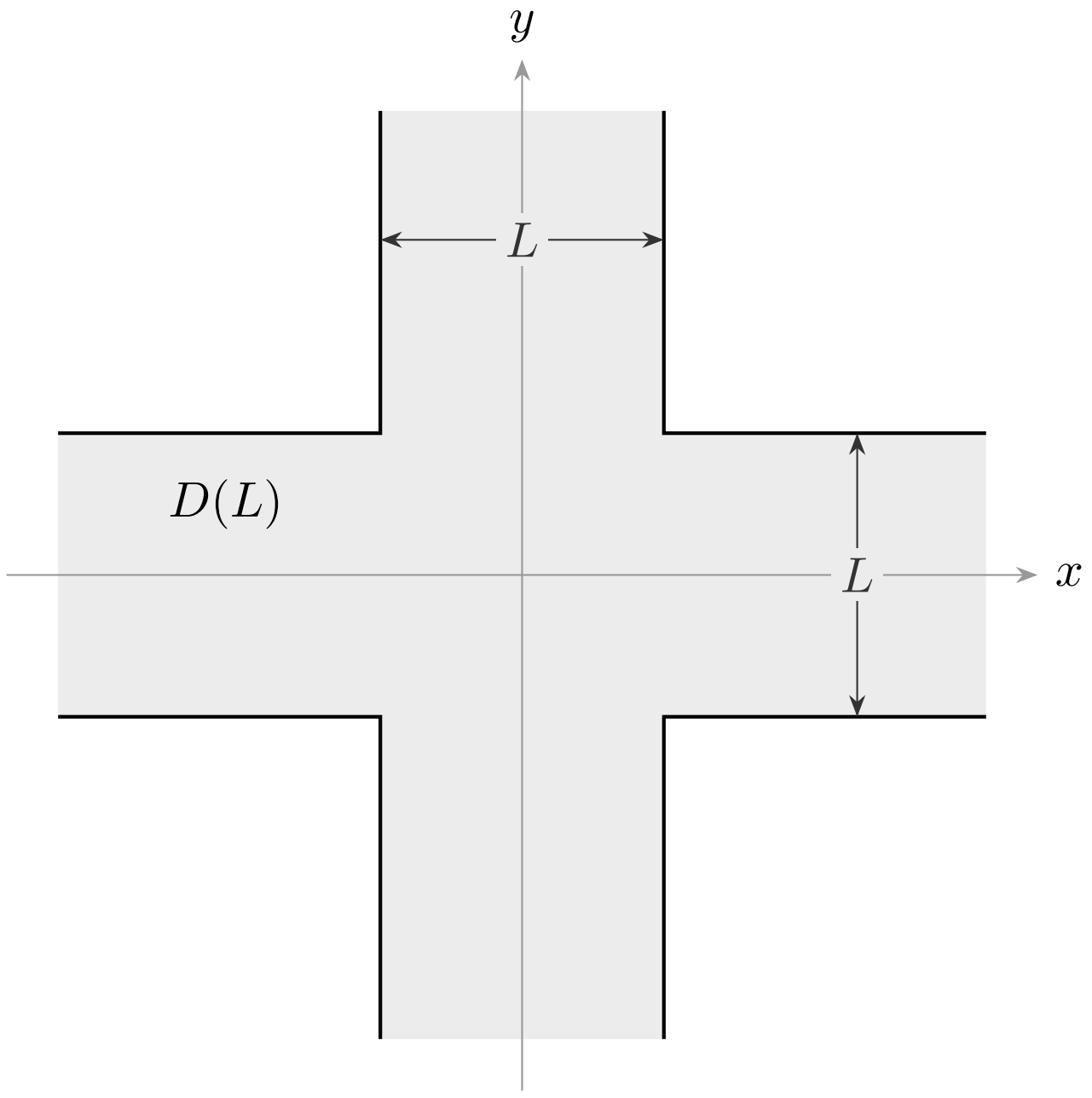
On $D(L)$, we investigate the existence and asymptotic profile of stationary solutions $u \colon D(L) \to \R$ satisfying the non-linear elliptic equation
\begin{equation}  \label{eqNlsCross}
    \begin{cases}
        -\Delta u + \omega u = \abs{u}^{p-2} u & \text{in } D(L), \\
        \dfrac{\partial u}{\partial \nu} = 0 & \text{on } \partial D(L),
    \end{cases}
\end{equation}
where $\omega > 0$ is a fixed frequency parameter, $p \in (2, \infty)$, and $\nu$ denotes the outward unit normal vector on the boundary $\partial D(L)$. We denote by $u_L$ solutions to \eqref{eqNlsCross} for some $L >0$.

Our main objective is to study \textit{the shrinking limit}
$L\to0^+$ of a family of positive stationary solutions $u_L$ to
\eqref{eqNlsCross}. The main aim is to show that the one-dimensional equation 
\[
    -\phi_e''+\omega\phi_e
    =|\phi_e|^{p-2}\phi_e,
    \qquad e=1,\ldots,4,
\]
governs their limiting profile on the four half-lines of a star graph, coupled by continuity and
the Kirchhoff condition at the vertex, defined for instance in \cite{zbMATH06132651}, see also Section~\ref{sec:1d_graph}. Because the two-dimensional domain $D(L)$ degenerates onto the underlying metric graph, with the arms shrinking only in the transverse direction and the central
junction shrinking in both directions, giving area factors \(L\) and
\(L^2\), the underlying functional spaces vary singularly with the thickness parameter. To compare solutions at different widths, we rescale the arms and the central junction separately.

Before addressing this limit, one must first construct stationary solutions to \eqref{eqNlsCross} for each fixed thickness. To capture the broadest possible range of non-linearity exponents, we choose to minimize the action functional 
\begin{equation*}
    S_{\omega,L}(u) = \frac{1}{2} \int_{D(L)} \left( \abs{\nabla u}^2 + \omega |u|^2 \right) \dxdy - \frac{1}{p} \int_{D(L)} \abs{u}^p \dxdy
\end{equation*}
constrained to nonzero functions satisfying $\langle S'_{\omega,L}(u),u \rangle =0$, that is on the Nehari manifold. Indeed, by the Euler-Lagrange theory, constrained critical points of the action on the  Nehari manifold satisfy \eqref{eqNlsCross}.  However, this strategy immediately introduces a severe functional-analytic obstacle. Because the cross-shaped domain is unbounded but lacks continuous translation invariance, minimizing sequences are susceptible to a loss of compactness. Indeed, the mass can drift infinitely far along the arms, and one cannot recover convergence by simply recentering the sequence as in the classical concentration-compactness scheme \cite{Lions84}.

This mechanism appears already on noncompact metric graphs as the \emph{runaway} alternative, as shown in
\cite[Lemma~3.3]{AdCaFiNo14} and
\cite[Section~3]{CaFiNo17}. And indeed, the limiting star graph, the graph made of two lines intersecting at zero, illustrates the consequences of
this phenomenon as for this graph, the infimum of the action over the
Nehari manifold is not attained and can be approached by soliton profiles
escaping along one edge \cite[Sections~3 and~5]{AdCaFiNo14Nehari}.
Nevertheless, stationary states do exist
\cite{AdCaFiNo12}, see also Section \ref{sec:1d_graph} below.

On the cross, as we show in Remark~\ref{rem:runaway_sequences} below, minimization exhibits an analogous family of escaping competitors. This motivates restricting to a cone of symmetric functions that are non-increasing along outward directions. The resulting compactness yields constrained minimizers; a further argument, based on the Neumann resolvent preserving the cone, shows that these minimizers solve  \eqref{eqNlsCross}. Similar methods have been used, for instance, in \cite{SeTi11, AlMo20, Moameni20}. In this setting, establishing the resolvent preservation introduces an elliptic regularity issue at the re-entrant corners, which we discuss below in Theorem \ref{thm:existence}.

Once this family of stationary solutions has been constructed, we study
its behavior as the width $L$ tends to zero. We prove strong convergence,
in a sense made precise in Theorem \ref{thm:main_shrinking_limit}, to the symmetric soliton on the limiting
four-edge star graph, with continuity and Kirchhoff conditions at the
central vertex defined in Section \ref{sec:1d_graph}.

\subsection{Presentation of the results}
We start by establishing the existence of stationary solutions for any $L >0$
\begin{theorem}
    \label{thm:existence}
    For all $L >0$ and $\omega >0$, there exist a solution $u_L \in H^1(D(L))$ to \eqref{eqNlsCross}. Moreover, $u_L\in H^2(D(L))$, it is symmetric with respect to axial symmetries $(x,y) \to (x,-y)$ and $(x,y) \to (-x,y)$, diagonal symmetry $(x,y) \to (y,x)$, and positive. There exists $M >0$ such that
    $$
0< u_L(x,y) \le M e^{-\sqrt{\omega/2} |(x,y)|}
$$
for $|(x,y)| \geq R(p, \omega) >0$. Finally, for $L\in(0,1]$ both $M$ and $R$ are independent of $L$.
\end{theorem}

As discussed above, a natural first approach to establish this existence result is to minimize the action $S_{\omega, L}$ on the unrestricted Nehari manifold. However, as the limiting star graph illustrates,
the infimum may fail to be attained because a minimizing sequence
escapes along a single arm. The cross admits analogous escaping
competitors (Remark~\ref{rem:runaway_sequences}), so a direct
minimization argument would have to rule out this possibility. 

To prevent escaping, we restrict the minimization problem to a closed cone of functions that are monotonically non-increasing along the outward directions. By enforcing this decay, we recover the compact embedding into $L^p$ space for any $p > 2$, guaranteeing the existence of a constrained minimizer $u_L$. 

However, because $u_L$ minimizes over a constrained cone rather than the full space $H^1(D(L))$, it initially satisfies only a variational inequality. To prove that $u_L$ is an unconstrained weak solution to the Euler-Lagrange equation \eqref{eqNlsCross}, we show that the linear Neumann resolvent operator preserves the cone. 

To establish the resolvent invariance, we enforce the axial symmetry $u(x,y) = u(|x|,|y|)$. This restricts the problem to the first-quadrant sub-domain, which we denote by $D_+(L)$; see Figure \ref{fig:domain_decomposition} below. A difficulty remains, however, in proving the resolvent invariance. In particular, to prove monotonicity preservation, one typically differentiates the weak equation associated with the resolvent problem (i.e., \eqref{eqNlsCross} with a forcing term on the right instead of the nonlinearity). One obtains an equation for $w = \partial_x u_L$ and tests it against its positive part $w^+$. This requires the resolvent solution to possess $H^2$ regularity so that $w \in H^1$ and the pairing is justified. 
But the domain $D_+(L)$ features a re-entrant corner with an internal angle of $3\pi/2$. As detailed in elliptic regularity theory for non-smooth domains (see, for example, Grisvard \cite[Section $8$]{Gr11}), the Neumann problem on a non-convex polygon generates a local singularity near this corner. Consequently, the resolvent operator suffers a degeneracy, mapping $L^2(D_+(L))$ only into the fractional Sobolev space $H^{5/3 - \epsilon}(D_+(L))$. In particular, $w$ need not belong to $H^1$, so the standard pairing with $w^+$ is mathematically unjustified. At the end of Section \ref{sec:ExisSol} we provide a counterexample based on the corner singularity which prevents resolvent preservation in this case.

To bypass this regularity bottleneck, we impose the diagonal symmetry $u(x,y) = u(y,x)$. This reduces the fundamental domain to the octant slice 
\begin{equation}
    \label{eqOctant}
    \Pi(L) \coloneqq \left\{ (x,y) : 0 < y < L/2, \ x > y \right\}.
\end{equation}
This diagonal bisection replaces the two corners $\pi/2$ and $3\pi/2$ with two new corners of $\pi/4$ and $3\pi/4$, making the domain convex. This eliminates the corner singularity, elevating the Neumann resolvent back to a standard $L^2 \to H^2$ mapping \cite[Theorem 3.2.1.3]{Gr11}.

With $H^2$ regularity, the differentiated resolvent identity is justified. The resolvent invariance property holds and implies that the constrained minimizer is a weak solution to \eqref{eqNlsCross}. Having established existence, we then study the uniform bounds and asymptotic properties of these solutions.

In the second part of the paper, we consider the shrinking limit for these constrained minimizers. Let $\mathcal{G}$ denote the underlying star graph made of four semi-infinite edges sharing the vertex identified with zero. On it, we define the profile 
\[
 \Phi_{\omega,\mathcal G}
 \coloneqq
 \bigl(\phi_\omega|_{\mathbb R^+},
       \phi_\omega|_{\mathbb R^+},
       \phi_\omega|_{\mathbb R^+},
       \phi_\omega|_{\mathbb R^+}\bigr)
\]
where $\phi_\omega$ is given by the one-dimensional positive soliton defined in \eqref{eqSechProf}.  Set $\Phi_{\omega,L}\coloneqq\mathcal C_L(\Phi_{\omega,\mathcal G})$ the cross extension of this profile, defined in \eqref{eqCrossExteGraph}. It is a construction that starts from a graph profile and trivially extends the edges in the transverse variable, while expanding the vertex value to the entire central junction. Let $s_{\omega}(L)$ denote the constrained minimization problem of the action on the Nehari manifold intersected with the cone in \eqref{eqCone}; see \eqref{eqMinProb}. In the same way, we denote by $s_{\omega, \mathcal{G}}$ the minimum level of the action on the Nehari manifold intersected with the set $K_{\mathrm{sym}}(\mathcal G)$ of symmetric and monotonically decreasing functions on the star graph; see Section \ref{sec1DProblem}. In particular, we show that 
$$s_{\omega, \mathcal{G}} = S_{\omega,\mathcal{G}}(\Phi_{\omega,\mathcal G})$$
holds; that is, this profile is a minimizer in the symmetric cone. Then we establish the following:

\begin{theorem}[Global strong convergence for the shrinking limit]
\label{thm:main_shrinking_limit}
Fix $\omega>0$ and $p>2$.  For $L>0$, let $u_L$ be a minimizer to $s_{\omega}(L)$ found in Theorem \ref{thmExistMinim}. 
Let $\Phi_{\omega,L}=\mathcal C_L(\Phi_{\omega,\mathcal G})$ the cross extension defined in \eqref{eqCrossExte}. Then, as $L\to0^+$,
\begin{equation*}
\lim_{L\to0}\frac{s_{\omega}(L)}{L}
=
s_{\omega,\mathcal G}
\end{equation*}
and
\[
\frac{1}{L}
\|u_L-\Phi_{\omega,L}\|_{H^2(D(L))}^2\to 0.
\]
\end{theorem}

The symmetries in the cone identify the four rescaled arms, so that the problem has a
single longitudinal average. The energy is of order $L$, whereas both
the transverse oscillation on an arm and the rescaled energy of the central
square vanish. The scale-sensitive Sobolev inequality in Lemma \ref{lem:thin_Sobolev} prevents the common
arm profile from disappearing in the limit. Its limit is therefore a
non-trivial symmetric monotone state on the star graph; the graph
variational characterization identifies it with the half-soliton on every
edge, which in particular satisfies the Kirchhoff condition. Finally, the
diagonal symmetry removes the leading re-entrant-corner singularity and
supplies a uniform $H^2$ estimate. This upgrades the energy convergence to $H^2$ convergence.

The cone constraint also plays a role in the shrinking limit. Suppose
that the unrestricted Nehari infimum were attained for every \(L>0\).
By Remark~\ref{rem:runaway_sequences}, the action of these minimizers,
divided by \(L\), would be at most the action of the line soliton.
Every positive stationary state on the four-edge star graph, however,
has twice that action. Thus even a family of unrestricted minimizers
could not converge in energy to a positive stationary state on the
graph: its profiles could drift along an arm as \(L\to0\).

\subsection{Discussion of the results and comparison with the literature}

We note that extending this analysis beyond intersections with $N > 2$ strips at a single junction to general unbounded networks is not trivial. Our proof of the resolvent preservation property relies essentially on the cross's reflection and diagonal symmetries and is not easily adaptable to a general network. On the other hand, without a suitable monotonicity constraint, the compactness mechanism  no
longer applies, and \textit{runaway} behavior must be addressed. The corresponding graph problems already show that global variational minimization is sensitive to graph topology. For energy minimization at prescribed mass in the subcritical range
\(2<p<6\), Adami, Serra and Tilli
\cite{AdSeTi15, AdSeTi16} identify a topological obstruction for
noncompact graphs. Under the assumption that every edge belongs to a trail joining two
distinct ends at infinity, the
infimum then equals the line-soliton energy and is not attained,
except on the real line and the special tower-of-bubbles
graphs. Outside this class, existence
cannot be decided from topology alone: edge lengths and the
prescribed mass may produce threshold phenomena
\cite{AdSeTi16}. At fixed frequency, analogous obstructions
and existence results for action minimization
on the Nehari manifold are established by De Coster et al. in
\cite{DeCoDoGaSeTr26}.

 These results concern global minimizers and do not rule out
other stationary states. They nevertheless indicate why a
generalization of the present variational construction is difficult for general thin domains, as one must find a mechanism that prevents loss of compactness while remaining compatible with the stationary equation.

The use of symmetry restrictions in variational problems is classically justified, under suitable hypotheses, by Palais' principle of symmetric criticality \cite{Palais79}. Monotonicity constraints require an additional argument to recover the full Euler-Lagrange equation. Variational constructions for Neumann problems with such constraints include the works of Serra and Tilli \cite{SeTi11} and Alves and Moameni \cite{AlMo20}. More generally, invariance conditions ensuring
that critical points on closed convex sets are unrestricted critical points have been studied by Marano and Mosconi in  \cite{MaMo15}, see also the references therein.
Moameni \cite{Moameni20} developed a general framework based on a pointwise invariance condition, which provides the abstract setting for the argument used here.

Our construction on the cross builds on unpublished notes by Moameni \cite{MoameniNotes} shared with the author, from which the monotonicity-cone approach and the compactness argument originate. Here, we address the re-entrant corners by incorporating diagonal symmetry into the cone, and formulate the minimization problem on the Nehari manifold to connect the existence argument with the subsequent analysis of the shrinking limit.

Dimensional reduction for nonlinear Schr\"odinger equations has been
investigated in several settings with a product structure. In the
mass-constrained framework, Terracini, Tzvetkov and Visciglia
\cite{TeTzVi14} study ground states on $\mathbb R^n\times M$, where
$M$ is a compact Riemannian manifold, and show that minimizers are independent of the compact
variable when the prescribed mass is sufficiently small. Recent new results in this direction are obtained by Pierotti, Verzini,  and Yu \cite{PiVeYu26}
 in the mass-supercritical regime and Soave, Verzini, and Villata \cite{SoVeVi26} for products of $\mathbb R^N$ and compact metric graphs.
For Neumann strips, Berestycki and Wei \cite{BeWe10} establish a transition between one-dimensional and transversely dependent least-energy solutions as the width varies. Le Coz and Shakarov have also studied dimensional reduction on fractured strips \cite{CoSh26Fr} and graph-based open books \cite{CoSh25Bk}, in the energy- and action-minimizing settings, respectively. Under the corresponding hypotheses, these results yield independence of the transverse variable for sufficiently small widths of energy and action minimizers.

A related phenomenon occurs for traveling waves with nonzero
conditions at infinity. For the Gross-Pitaevskii equation on
$\mathbb R\times\mathbb T$ with nonzero background, de Laire, Gravejat and Smets
\cite{deGrSm24} and   Mari\c{s} and Mur
\cite{MaMu24}  obtain one-dimensional minimizing profiles when
the transverse period is sufficiently small; see also Aftalion, Gravejat and Sandier \cite{MR5074438} and Aftalion and Nguyen \cite{zbMATH08197506} for recent related results.

An important feature of product geometries is the separation
between longitudinal and transverse variables. Rescaling the
compact factor transfers the dependence on its size to the
coefficients of a problem on a fixed space. This structure allows
variational comparisons and rigidity arguments to establish convergence and transverse independence below a
certain threshold. The cross-shaped domain does not admit such a global transverse rescaling. Although the isotropic dilation $D(L)=L D(1)$ fixes the domain, it also rescales the longitudinal variables and therefore does not preserve the limiting soliton's scale. We instead rescale the arms and the central
junction separately and analyze their coupling. Our result thus differs from those above as it is an
asymptotic reduction to a graph profile as $L\to0$, rather than an
identification of finite-width solutions with a graph state
extended constantly in the transverse direction.

For domains with junctions, Kosugi \cite{Ko00, Ko02} provides a closer geometric precedent, studying semilinear elliptic equations on thin networks. However, these works assume a bounded domain, smooth lateral boundaries, and dissipativity conditions on
the nonlinearity at large amplitudes, which supply uniform pointwise bounds by the maximum principle. Moreover, they do not consider a variational approach. 

In contrast, we deal with unbounded arms, a focusing power nonlinearity, and unsmoothed re-entrant corners. In particular, the uniform
pointwise bounds needed for the limiting analysis do not follow from the same argument. Moreover, to study the shrinking limit, we combine the convergence of the variational levels to upgrade weak $H^1$ convergence to strong convergence in the suitably rescaled energy spaces and in $H^2$.

\subsection{Outline of the paper}
The remainder of the paper is organized as follows. In Section~\ref{secPrelimin},
we recall the variational and functional framework on the four-edge
star graph. Section \ref{sec:ExisSol} is devoted to the existence theory on the
cross-shaped domain: we prove the existence of a minimizer in
Theorem~\ref{thmExistMinim}, establish in
Lemma~\ref{lem:cone_invariance} that this minimizer is a 
solution to \eqref{eqNlsCross}, and discuss the role of the cone constraints, including
a counterexample to resolvent invariance without diagonal symmetry. The corresponding
one-dimensional variational problem on the star graph is treated in
Subsection~\ref{sec1DProblem}. In Section~\ref{sec:properties}, we
establish qualitative properties of the stationary states constructed
above and prove Theorem~\ref{thm:existence}. Finally, in
Section~\ref{sec:shrinking_limit}, we study the limit as the width of
the cross tends to zero and prove Theorem~\ref{thm:main_shrinking_limit}. For the rest of this work, we will refer to $D(L)$ as \textit{the cross}.

\section{Preliminaries}\label{secPrelimin}

\subsection{The one-dimensional graph problem}
\label{sec:1d_graph}

To prepare for the shrinking limit as $L \to 0$, we introduce the limiting one-dimensional problem on a star graph $\mathcal{G}$. We define the star graph $\mathcal{G}$ as the union of four infinite half-lines (edges) $e_i \simeq [0, \infty)$ for $i \in \{1, 2, 3, 4\}$, joined at a single central vertex $\nu = 0$.

A function $v$ on $\mathcal{G}$ can be represented as a set $v = (v_1, v_2, v_3, v_4)$ where each $v_i \colon [0, \infty) \to \R$. The Sobolev space $H^1(\mathcal{G})$ is defined as the set of functions $v \in \bigoplus_{i=1}^4 H^1(e_i)$ satisfying the continuity condition at the central vertex:
\begin{equation}
    H^1(\mathcal{G}) \coloneqq \{v \in \bigoplus_{i=1}^4 H^1(e_i): v_1(0) = v_2(0) = v_3(0) = v_4(0)\}.
    \label{eqGraphConti}
\end{equation}
The corresponding one-dimensional action functional $S_{\omega, \mathcal{G}} \colon H^1(\mathcal{G}) \to \R$ is given by
\begin{equation*}
    S_{\omega, \mathcal{G}}(v) \coloneqq \sum_{i=1}^4 \left( \frac{1}{2} \int_0^\infty \left( \abs{v_i'}^2 + \omega |v_i|^2 \right) dx - \frac{1}{p} \int_0^\infty \abs{v_i}^p dx \right),
\end{equation*}
Critical points of $S_{\omega, \mathcal{G}}$ correspond to weak solutions of the limiting stationary nonlinear Schr\"odinger equation on $\mathcal{G}$. In particular, each $v_i$ satisfies
\begin{equation}
    \label{eq1DNLS}
    -\partial_x^2 v_i + \omega v_i = |v_i|^{p-2} v_i ,
\end{equation}
and at the vertex $\nu = 0$, it satisfies the Kirchhoff boundary condition:
\begin{equation}
    \sum_{i=1}^4 v_i'(0) = 0,
    \label{eqKirchhoff}
\end{equation}
where $v_i'(0)$ denotes the outgoing derivative from the vertex along the edge $e_i$.

Outside of the trivial solution, the set of solutions to \eqref{eq1DNLS} in $H^1(\mathcal G)$ satisfying \eqref{eqKirchhoff} is constructed by attaching the one-dimensional half-soliton profile to each of the four edges. Let $\phi_\omega$ be the unique positive even ground state solution to the one-dimensional problem 
$$-\phi'' + \omega \phi = \abs{\phi}^{p-2}\phi$$
on $\R$, given explicitly by
\begin{equation}    \label{eqSechProf}
    \phi_\omega(x) \coloneqq \left( \frac{p \omega}{2} \right)^{\frac{1}{p-2}} \sech^{\frac{2}{p-2}} \left( \frac{(p-2)\sqrt{\omega}}{2} x \right).
\end{equation}
The only positive solutions on the star graph $\mathcal{G}$ are then defined by 
\begin{equation}
   \phi_{\omega,\mathcal{G}} := (\phi_\omega(x - a_i))_i \quad \text{for } x \in [0, \infty) \text{ and } i \in \{1, 2, 3, 4\},
    \label{eqGraphGS}
\end{equation}
for some parameters $a_i \in \R$, such that $|a_i| = a$ for any $i$, which implies the continuity at the vertex \eqref{eqGraphConti} and (up to a permutation of the indices) $a_1 = a_2 = - a_3 = -a_4$, which implies \eqref{eqKirchhoff}.

\section{Existence of stationary solutions}\label{sec:ExisSol}

\subsection{Definitions of sub-domains}
\label{sec:variational}
We define the central square $q(L)$ and the four disjoint semi-infinite strips corresponding to the arms of the cross $D(L)$ by 
\begin{align}
   \nonumber q(L) &\coloneqq (-L/2, L/2) \times (-L/2, L/2), \\ \label{eqH+}
    h^+(L) &\coloneqq (L/2, \infty) \times (-L/2, L/2), \quad &h^-(L) &\coloneqq (-\infty, -L/2) \times (-L/2, L/2), \\ \nonumber
    v^+(L) &\coloneqq (-L/2, L/2) \times (L/2, \infty), \quad &v^-(L) &\coloneqq (-L/2, L/2) \times (-\infty, -L/2),\\ 
    h(L) &\coloneqq \R \times (-L/2, L/2) , \quad & v(L) &\coloneqq  (-L/2, L/2) \times \R. \label{eqHandV}
\end{align}
To properly enforce symmetry, it is convenient to introduce the first-quadrant portion of the cross :
\begin{equation}
    D_+(L) \coloneqq \left( (0, \infty) \times \left(0, \frac{L}{2}\right) \right) \cup \left( \left(0, \frac{L}{2}\right) \times (0, \infty) \right).
    \label{eqFirstQuad}
\end{equation}

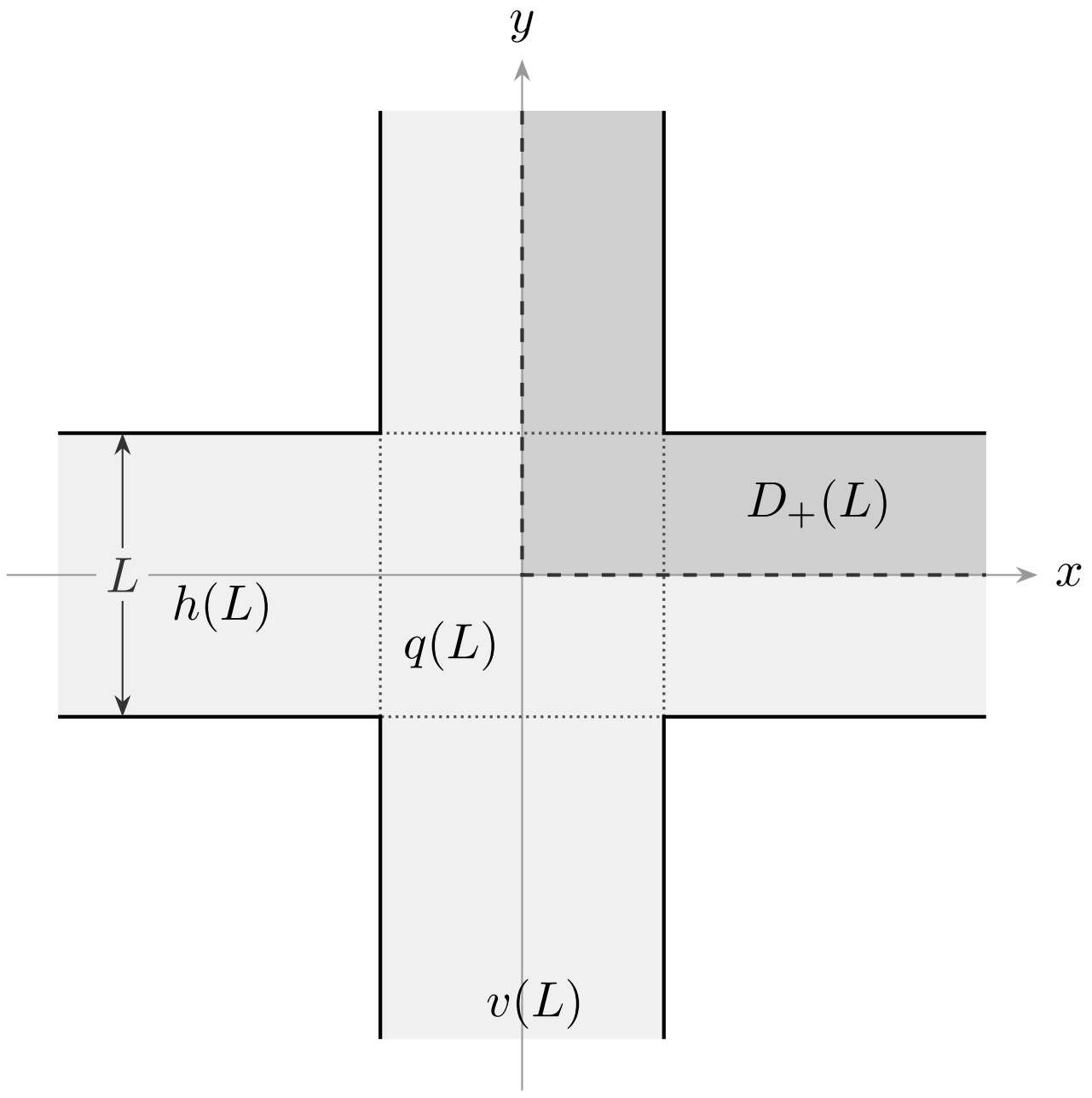
\begin{figure}[htbp]
    \centering
    \begin{tikzpicture}[
        scale=1.0,
        >=Stealth,
        boundary/.style={thick, black},
        symboundary/.style={thick, black!80, dashed},
        interface/.style={semithick, black!70, densely dotted},
        axis/.style={->, thin, gray!80},
        dim/.style={<->, thin, black!80}
    ]
        \def\a{1.1}
        \def\m{3.6}

        \fill[gray!12] 
            (-\m, -\a) -- (-\a, -\a) -- (-\a, -\m) -- 
            ( \a, -\m) -- ( \a, -\a) -- ( \m, -\a) -- 
            ( \m,  \a) -- ( \a,  \a) -- ( \a,  \m) -- 
            (-\a,  \m) -- (-\a,  \a) -- (-\m,  \a) -- cycle;

        \fill[gray!38] 
            (0, 0) -- (\m, 0) -- (\m, \a) -- 
            (\a, \a) -- (\a, \m) -- (0, \m) -- cycle;

        \draw[axis] (-\m - 0.4, 0) -- (\m + 0.4, 0) node[right, black] {$x$};
        \draw[axis] (0, -\m - 0.4) -- (0, \m + 0.4) node[above, black] {$y$};

        \draw[interface] (-\a, -\a) rectangle (\a, \a);

        \draw[symboundary] (0, \m) -- (0, 0) -- (\m, 0);

        \draw[boundary] (-\m,  \a) -- (-\a,  \a) -- (-\a,  \m);
        \draw[boundary] ( \m,  \a) -- ( \a,  \a) -- ( \a,  \m);
        \draw[boundary] (-\m, -\a) -- (-\a, -\a) -- (-\a, -\m);
        \draw[boundary] ( \m, -\a) -- ( \a, -\a) -- ( \a, -\m);

        \draw[dim] (-3.1, -\a) -- (-3.1, \a) node[midway, fill=gray!12, inner sep=2pt] {$L$};

        \node at ( 2.3,  0.55) {$D_+(L)$};
        \node at (-0.55, -0.55) {$q(L)$};
   
        \node[left]  at (-1.8, -0.25) {$h(L)$};
        \node[below] at (0.1, -3)  {$v(L)$};
    \end{tikzpicture}
    \caption{Decomposition of $D(L)$ into the central square $q(L)$, the horizontal and vertical arms $h(L)$ and $v(L)$, and the first-quadrant subdomain $D_+(L)$.}
    \label{fig:domain_decomposition}
\end{figure}
   
We then introduce the functional space $\mathcal{K}(L) \subset H^1(D(L))$ defined by
\begin{definition}[Symmetric monotonic cone]
Let \(\tau(x,y)\coloneqq(y,x)\). We set
\begin{equation}
    \label{eqCone}
    \mathcal K(L)\coloneqq
\left\{
u\in H^1(D(L)):
\begin{array}{l}
u(-x,y)=u(x,y)=u(x,-y),\\
u\circ\tau=u,\\
\partial_xu\leq0,\quad \partial_yu\leq0
\quad\text{a.e.\ in }D_+(L)
\end{array}
\right\}.
\end{equation}
\end{definition}
Functions in $\mathcal{K}(L)$ are even in both variables, symmetric about the diagonal, and non-increasing away from the origin along the unbounded directions. Moreover, any $u\in \mathcal{K} (L)$ is nonnegative since it is globally in $H^1$.  

\begin{definition}[Cross extension]
\label{def:cross_extension}
For any $L > 0$, we define the \emph{cross extension} operator $\mathcal{C}_L: C^0([0, \infty)) \to C^0(D(L))$ mapping a 1D profile $\phi$ to $\phi_L \coloneqq \mathcal{C}_L(\phi)$ as the continuous piecewise function:
\begin{equation}
    \label{eqCrossExte}
    \phi_L(x,y) \coloneqq \begin{cases}
        \phi(0) & \text{for } (x,y) \in q(L), \\
        \phi(|x| - L/2) & \text{for } (x,y) \in h^\pm(L), \\
        \phi(|y| - L/2) & \text{for } (x,y) \in v^\pm(L).
    \end{cases}
\end{equation}
Similarly, for a star graph function $v \in H^1(\mathcal{G})$, its cross extension is defined as
\begin{equation}
    \label{eqCrossExteGraph}
    v_L(x,y) \coloneqq \begin{cases}
        v(0) & \text{for } (x,y) \in q(L), \\
        v_1(|x| - L/2) & \text{for } (x,y) \in h^+(L), \\
        v_2(|y| - L/2) & \text{for } (x,y) \in v^+(L), \\
        v_3(|x| - L/2) & \text{for } (x,y) \in h^-(L), \\
        v_4(|y| - L/2) & \text{for } (x,y) \in v^-(L).
    \end{cases}
\end{equation}
\end{definition}

\subsection{The Nehari manifold and minimization}

The action functional associated with \eqref{eqNlsCross} on $D(L)$ is defined for $u \in H^1(D(L))$ by $S_{\omega,L} \colon  H^1(D(L)) \to \R$:
\begin{equation}\label{eqActionFun}
    S_{\omega,L}(u) \coloneqq \frac{1}{2} \int_{D(L)} \left( \abs{\nabla u}^2 + \omega |u|^2 \right) \dxdy - \frac{1}{p} \int_{D(L)} \abs{u}^p \dxdy.
\end{equation}
The minimization problem is constrained to the Nehari manifold restricted to $\mathcal{K}(L)$
\begin{equation*}
    \mathcal{N}(L) \coloneqq \left\{ u \in \mathcal{K}(L) \setminus \{0\} : I_{\omega,L}(u) = 0 \right\},
\end{equation*}
where
\[
I_{\omega,L}(u) \coloneqq \langle S_{\omega,L}'(u), u \rangle = \int_{D(L)} (\abs{\nabla u}^2 + \omega |u|^2) \dxdy - \int_{D(L)} \abs{u}^p \dxdy.
\]

We define the minimal \textit{constrained} ground state energy level $s_\omega(L)$ as the minimization problem:
\begin{equation}\label{eqMinProb}
    s_\omega(L) \coloneqq \inf_{u \in \mathcal{N}(L)} S_{\omega,L}(u).
\end{equation}

\begin{lemma}\label{lemInfDef}
For every $u \in \mathcal{N}(L)$, the following holds:
\begin{equation}
    S_{\omega,L}(u) = \frac{p-2}{2p} \int_{D(L)} \left( \abs{\nabla u}^2 + \omega |u|^2 \right) \dxdy = \frac{p-2}{2p} \int_{D(L)} \abs{u}^p \dxdy.
    \label{eqActionPotential}
\end{equation}
Furthermore, the infimum level $s_\omega(L)$ verifies
\begin{equation}\label{eqMinProbPot}
    s_\omega(L) = \frac{p-2}{2p} \inf \left\{ \| u \|_{L^p(D(L))}^p : u \in \mathcal{K}(L) \setminus \{0\}, \, I_{\omega,L}(u) \le 0 \right\}.
\end{equation}
\end{lemma}

\begin{proof}
Identity \eqref{eqActionPotential} follows from the fact that on $\mathcal{N}(L)$,
\begin{equation*}
    S_{\omega,L}(u) - \frac{1}{2} I_{\omega,L}(u) = \frac{p-2}{2p} \int_{D(L)} \abs{u}^p \dxdy.
\end{equation*}
To verify \eqref{eqMinProbPot}, note that 
$$\mathcal{N}(L) \subset \{ u \in \mathcal{K}(L) \setminus \{0\} : I_{\omega,L}(u) \le 0 \},$$
therefore
\[
\inf_{I_{\omega,L} \le 0} \frac{p-2}{2p}\|u\|_{L^p(D(L))}^p \le s_\omega(L).
\]
Conversely, let $u \in \mathcal{K}(L) \setminus \{0\}$ satisfy $I_{\omega,L}(u) < 0$. Define the Nehari projection scalar
\begin{equation}
    \label{eqNehariProj}
    \pi(u) \coloneqq \left( \frac{\int_{D(L)} (\abs{\nabla u}^2 + \omega |u|^2) \dxdy}{\int_{D(L)} \abs{u}^p \dxdy} \right)^{\frac{1}{p-2}},
\end{equation}
and observe that $\pi(u) u \in \mathcal{N}(L)$. Since $I_{\omega,L}(u) < 0$, it follows that $\pi(u) \in (0,1)$. Consequently,
\[
s_\omega(L) \le \frac{p-2}{2p}\| \pi(u) u\|_{L^p(D(L))}^p < \frac{p-2}{2p}\|u\|_{L^p(D(L))}^p,
\]
which shows that minimizing over $I_{\omega,L}(u) \le 0$ yields the exact same infimum level $s_\omega(L)$.
\end{proof}

\subsection{Compactness and existence of solutions on the cross}

To rigorously establish compactness on the full symmetric domain, we first isolate the behavior on the first quadrant. We define the restricted function space on the first quadrant $D_+(L)$ defined in \eqref{eqFirstQuad}.
\begin{definition}
We set
\begin{equation*}
    \mathcal K_+(L)\coloneqq
\left\{
u\in H^1(D_+(L)):
\begin{array}{l}
u\circ\tau=u,\\
\partial_xu\leq0,\quad \partial_yu\leq0
\quad\text{a.e. in }  D_+(L)
\end{array}
\right\}.
\end{equation*}
\end{definition}

By symmetry, one can extend functions in $\mathcal K_+(L)$ to those in $\mathcal K(L)$ and, conversely, restrict functions from $\mathcal K(L)$ to $\mathcal K_+(L)$, as remarked in the following Lemma.

\begin{lemma}[Symmetric extension]\label{lemProjection}
The extension operator 
\[
\mathcal{E} : H^1(D_+(L))\longrightarrow H^1(D(L)),  \quad (\mathcal{E}u)(x,y) \coloneqq u(|x|,|y|)
\]
is well-defined and $\mathcal K_+(L) \to \mathcal{K}(L)$. Conversely, every $w \in \mathcal{K}(L)$ is of the form $w(x,y) = u(|x|,|y|)$ for a unique $u \in \mathcal{K}_+(L)$, namely $u = w|_{D_+(L)}$. Moreover, it holds
\[
\|w|_{D_+(L)}\|_{H^1(D_+(L))} \le \|w\|_{H^1(D(L))},
\]
and
\[
\|\mathcal{E}u\|_{L^q(D(L))}^q \le 4 \|u\|_{L^q(D_+(L))}^q \quad \text{for every } q \ge 1.
\]
\end{lemma}

\begin{proof}
The restriction statement is immediate. Conversely, if $u \in \mathcal{K}_+(L)$, define $w(x,y) \coloneqq u(|x|,|y|)$. Clearly, one has $w \in \mathcal{K}(L)$. Moreover, one directly gets
$$\|w|_{D_+(L)}\|_{H^1(D_+(L))} \le \|w\|_{H^1(D(L))}.$$
and by symmetry
\[
\|w\|_{L^q(D(L))}^q \le 4 \|u\|_{L^q(D_+(L))}^q.
\]
\end{proof}

\begin{lemma}[Compactness on the first quadrant]\label{lemCompc1/4}
For every $q > 2$, the embedding
\[
\mathcal{K}_+(L) \hookrightarrow L^q(D_+(L))
\]
is compact.
\end{lemma}

\begin{proof}
Let $(u_n) \subset \mathcal{K}_+(L)$ be bounded in $H^1(D_+(L))$ with $\|u_n\|_{H^1(D_+(L))} \le M$ for all $n$. We prove that $(u_n)$ admits a strongly convergent subsequence in $L^q(D_+(L))$.

We start by showing a uniform tail estimate on the horizontal part of $D_+(L)$, that is, $[0, \infty) \times [0, L/2]$. For each $n$, define
\[
F_n(x) := \left\| |u_n(x,y)| \right\|_{L^\infty_y(0,L/2)}^2, \quad x \ge 0.
\]
Since $\partial_x u_n \le 0$ a.e. and $u_n(\cdot, y) \in L^2(0,\infty)$, we   have $u_n \ge 0$ a.e. (otherwise $u_n(x,y) \le c < 0$ for large $x$, contradicting square integrability). Hence $x \mapsto u_n(x,y)$ is non-negative and non-increasing, so $F_n$ is non-increasing on $[0,\infty)$.

Fix $x \ge 0$. By the mean value theorem for integrals, there exists $y_n \in [0, L/2]$ such that
\begin{equation}\label{eqAverage}
    |u_n(x,y_n)|^2 \leq \frac{2}{L} \int_0^{L/2} |u_n(x,y)|^2 dy. 
\end{equation}
For any $y \in [0, L/2]$, the Cauchy-Schwarz inequality yields
\begin{align}
   |u_n(x,y)|^2 - |u_n(x,y_n)|^2 &\leq \int_{y_n}^y \partial_z |u_n(x,z)|^2  dz \leq 2  \int_{0}^{L/2} |u_n(x,z)| |\partial_z u_n(x,z)| dz \nonumber \\
   &\leq \int_{0}^{L/2} |u_n(x,z)|^2 dz + \int_{0}^{L/2} |\partial_z u_n(x,z)|^2 dz. \label{eq1DSobol}
\end{align}
Combining \eqref{eqAverage} and \eqref{eq1DSobol}, we obtain
\[
|u_n(x,y)|^2 \leq    \left(1+\frac2L \right)
\int_0^{L/2}|u_n(x,z)|^2\,dz + \int_0^{L/2}|\partial_yu_n(x,z)|^2\,dz,
\] 
Taking the supremum over $y$ to get $F_n(x)$ and integrating over $x \in [0,\infty)$ yields
\[
\int_{0}^{\infty} F_n(x) \, dx \le  \left(1+\frac2L \right)  \int_{[0, \infty) \times [0, L/2]} \left( |u_n|^2 + |\partial_y u_n|^2 \right) dx dy \le  \left(1+\frac2L \right) M^2.
\]
Since $F_n \ge 0$ is non-increasing, integration by parts guarantees 
$$x F_n(x) \le \int_{0}^{x} F_n(t) \, dt \le  \left(1+\frac2L \right)  M^2,$$
so 
\begin{equation}
    \label{eqPolyDecay}
    F_n(x) \le \left(1+\frac2L \right) \frac{  M^2}{x}.
\end{equation}
For any $q > 2$ and $R > 0$,
\[
\int_{R}^\infty \int_0^{L/2} |u_n|^q \, dy dx \le \int_{R}^{\infty} \int_{0}^{L/2} F_n(x)^{q/2} \, dy dx \lesssim  M^q \int_{R}^{\infty} x^{-q/2} \, dx.
\]
Since $q/2 > 1$, the right-hand side vanishes as $R \to \infty$, giving uniform tail decay on the horizontal section.

Applying a symmetric argument to the vertical section $[0, L/2] \times [0, \infty)$ with respect to the $y$-variable yields the same uniform tail decay. Setting $D_+(L)_R \coloneqq D_+(L) \cap ([0,R] \times [0,R])$, combining both sections gives
\[
\lim_{R \to \infty} \sup_n \int_{D_+(L) \setminus D_+(L)_R} |u_n|^q \, dx dy = 0.
\]

This yields strong $L^q$ convergence. Indeed, for fixed $R > 0$, the domain $D_+(L)_R$ is bounded and Lipschitz. As $(u_n)$ is bounded in $H^1(D_+(L)_R)$, the Rellich--Kondrachov theorem ensures that, up to a subsequence, $u_n \to u$ strongly in $L^q(D_+(L)_R)$. For any $\varepsilon > 0$, choosing $R$ sufficiently large controls the tails, and the local strong convergence on the compact core guarantees global convergence: $u_n \to u$ strongly in $L^q(D_+(L))$.
\end{proof}

\begin{proposition}[Compactness on the full cross]\label{prop:compactness_HL}
For every $q > 2$, the embedding
\[
\mathcal{K}(L) \hookrightarrow L^q(D(L))
\]
is compact.
\end{proposition}

\begin{proof}
Let $(w_n) \subset \mathcal{K}(L)$ be a sequence bounded in $H^1(D(L))$. For each $n$, set $u_n \coloneqq w_n|_{D_+(L)}$. By definition of $\mathcal{K}(L)$, each $u_n$ belongs to $\mathcal{K}_+(L)$, and by the restriction estimate in Lemma \ref{lemProjection}, $\|u_n\|_{H^1(D_+(L))} \le C_1 \|w_n\|_{H^1(D(L))} \le M$. 

By Lemma \ref{lemCompc1/4}, after passing to a subsequence, there exists $u \in L^q(D_+(L))$ such that $u_n \to u$ strongly in $L^q(D_+(L))$. Setting $w \coloneqq \mathcal{E}u$, the extension bound from Lemma \ref{lemProjection} applied to the difference $w_n - w = \mathcal{E}(u_n - u)$ yields
\[
\|w_n - w\|_{L^q(D(L))} \le C_2 \|u_n - u\|_{L^q(D_+(L))} \xrightarrow{n \to \infty} 0.
\]
Hence, $w_n \to w$ strongly in $L^q(D(L))$.
\end{proof}
Next we need the following closeness lemma. 

\begin{lemma}[Closeness of the cone]
\label{lem:cone_closed}
The cone $\mathcal K(L)$ is strongly and weakly closed in $H^1(D(L))$. 
\end{lemma}

\begin{proof}
Let $u_n\in\mathcal K(L)$ and suppose that
$u_n\to u$ strongly in $H^1(D(L))$. Since $u_n\geq0$, we have $u\geq0$ almost everywhere.
Moreover, let $g$ be any of the two coordinate reflections or the diagonal
reflection $\tau$. Composition with $g$ is an isometry of
$H^1(D(L))$. Since $u_n\circ g=u_n$ for every $n$, passing to the
limit gives
\[
u\circ g=u.
\]
Finally, for $j\in\{x,y\}$,  since $\partial_j u_n\leq0$ almost everywhere in $D_+(L)$, so is $\partial_j u\leq0$ almost everywhere in $D_+(L)$.
Hence $u\in\mathcal K(L)$. So $\mathcal K(L)$ is strongly closed. Since $\mathcal K(L)$ is convex, it is also weakly closed in $H^1(D(L))$.
\end{proof}

We now show existence of minimizers to \eqref{eqMinProb}.

\begin{theorem}[Existence of a constrained minimizer]
\label{thmExistMinim}
For every fixed $L > 0$ and $\omega > 0$, the minimization problem \eqref{eqMinProb} admits a non-trivial and nonnegative minimizer $u_L \in \mathcal{N}(L)$. 
\end{theorem}

\begin{proof}
Let $(u_n)_{n \in \N} \subset \mathcal{N}(L)$ be a minimizing sequence for $S_{\omega,L}$. We can assume $I_{\omega,L}(u_n) = 0$ for every $n$. Define 
$$
\|u_n\|_{H^1_\omega(D(L))}^2 \coloneqq \int_{D(L)} (|\nabla u_n|^2 + \omega |u_n|^2) \, dxdy.
$$
We have by \eqref{eqMinProbPot}
\[
\|u_n\|_{H^1_\omega(D(L))}^2  =   \frac{2p}{p-2} s_\omega(L) + o(1).
\]
This implies that $(u_n)$ is bounded in $H^1(D(L))$. Moreover, by the Sobolev embedding, we have 
$$\|u_n\|_{L^p(D(L))}^p \le C \|u_n\|_{H^1_\omega(D(L))}^p$$
for $p > 2$. Using the Nehari identity, this yields
\[
\|u_n\|_{H^1_\omega(D(L))}^2 = \|u_n\|_{L^p(D(L))}^p \le C \|u_n\|_{H^1_\omega(D(L))}^p.
\]
Since $u_n \neq 0$, dividing by $\|u_n\|_{H^1_\omega(D(L))}^2$ yields 
$$\|u_n\|_{L^p(D(L))}^p \ge C^{-2/(p-2)} > 0$$
for all $n$. Consequently, $s_\omega(L) > 0$. 
Since $(u_n)$ is bounded in $\mathcal{K}(L)$, the compact embedding $\mathcal{K}(L) \hookrightarrow L^p(D(L))$ from Proposition \ref{prop:compactness_HL} guarantees that, up to a subsequence, there exists $u_L \in \mathcal{K}(L)$ such that $u_n \rightharpoonup u_L$ weakly in $H^1(D(L))$ and $u_n \to u_L$ strongly in $L^p(D(L))$. By strong convergence,
\[
\|u_L\|_{L^p(D(L))}^p = \lim_{n \to \infty} \|u_n\|_{L^p(D(L))}^p = \frac{2p}{p-2} s_\omega(L) > 0,
\]
so $u_L \neq 0$. By weak lower semicontinuity of the $H^1_\omega$-norm,
\[
\|u_L\|_{H^1_\omega(D(L))}^2 \le \liminf_{n \to \infty} \|u_n\|_{H^1_\omega(D(L))}^2 = \lim_{n \to \infty} \|u_n\|_{L^p(D(L))}^p = \|u_L\|_{L^p(D(L))}^p.
\]
Hence, $I_{\omega,L}(u_L) \le 0$. By Lemma \ref{lemInfDef}, we get
$$\frac{p-2}{2p} \|u_L\|_{L^p(D(L))}^p = s_\omega(L).$$ 
As in the proof of Lemma \ref{lemInfDef}, this implies that $I_{\omega,L}(u_L) = 0$, that is 
$$\|u_L\|_{H^1_\omega(D(L))}^2 = \|u_L\|_{L^p(D(L))}^p.$$
By Lemma \ref{lem:cone_closed}, $u_L \in \mathcal{K}(L)$. Thus $u_L \in \mathcal{N}(L)$ and $S_{\omega,L}(u_L) = s_\omega(L)$. Finally, the non-negativity follows from the definition of $\mathcal{K}(L)$.
\end{proof}

\begin{remark}\label{rem:StricPos}
    Once Theorem \ref{thm:existence} is shown, one can conclude that $u_L$ is a solution to the elliptic problem \eqref{eqNlsCross}. Then, by interior regularity and the strong maximum principle, we conclude that $u_L$ in Proposition \ref{thmExistMinim} is strictly positive in $D(L)$.
\end{remark}

\subsection{Monotonicity and symmetry preservation of the resolvent} \label{sec:resolvent}

To complete the proof of Theorem \ref{thm:existence}, we demonstrate that the constrained minimizer $u_L$ is a true weak solution. We establish this by showing the resolvent operator preserves the symmetric monotonicity cone $\mathcal{K}(L)$. The general principle underlying this argument appears in \cite[Theorem~1.1]{Moameni20}.   The approach is based on showing that the Neumann resolvent  $\mathcal{R}_{\omega,L} \coloneqq (-\Delta+\omega)^{-1}$ preserves the cone $\mathcal{K}(L)$. In particular, we will show the following
\begin{equation}
\label{eqResolventInv}
u \in \mathcal{K}(L), \qquad \implies \qquad \mathcal{R}_{\omega,L}(u^{p-1}) \in \mathcal{K}(L).
\end{equation}
This property is sufficient to claim that the \textit{constrained} minimizer $u_L$ is actually a weak solution to the Euler-Lagrange equation associated with the action. The next proposition shows this.

\begin{proposition}
\label{prop:natural_cone}
Assume \eqref{eqResolventInv} holds. Let $u_L\in \mathcal{K}(L)$ be a nonnegative minimizer defined in Theorem \ref{thmExistMinim}. Then $u_L$ is a weak solution of \eqref{eqNlsCross} over the entire space $H^1(D(L))$; namely,
\[
\int_{D(L)} \left( \nabla u_L\cdot\nabla\psi+\omega u_L\psi-u_L^{p-1}\psi \right) dx\,dy = 0 \qquad \text{for every } \psi\in H^1(D(L)).
\]
\end{proposition}

\begin{proof}
Let
$$a_{\omega,L}(u,v) \coloneqq \int_{D(L)} \left(\nabla u\cdot\nabla v+\omega uv\right) dx\,dy.$$
For $\pi$ defined in \eqref{eqNehariProj}, let 
$$F(z) \coloneqq S_{\omega,L}(\pi(z)z)$$
for $z\in H^1(D(L))\setminus\{0\}$.   Since $u_L\in\mathcal{N}(L)$, we have $\pi(u_L)=1$ and $\langle S_{\omega,L}'(u_L),u_L\rangle = 0$. By the chain rule, for every $h\in H^1(D(L))$:
\begin{align*}
F'(u_L)[h] &= S_{\omega,L}'(u_L) \left[ h + \pi'(u_L)[h]u_L \right] \\
&= S_{\omega,L}'(u_L)[h] + \pi'(u_L)[h] \langle S_{\omega,L}'(u_L),u_L \rangle = S_{\omega,L}'(u_L)[h].
\end{align*}
Because $\mathcal{K}(L)$ is convex, the segment $u_L+\varepsilon(v-u_L)$ lies in $\mathcal{K}(L)$ for every $v\in \mathcal{K}(L)$ and $\varepsilon\in[0,1]$. The minimality of $u_L$ for $F$ on this convex set implies the variational inequality $S_{\omega,L}'(u_L)[v-u_L]\geq0$ for every $v\in \mathcal{K}(L)$.

We introduce the resolvent function $v_L \coloneqq \mathcal{R}_{\omega,L}(u_L^{p-1})$. By Lemma \ref{lem:cone_invariance}, $v_L\in \mathcal{K}(L)$. The resolvent satisfies 
$$a_{\omega,L}(v_L,\psi) = \int_{D(L)} u_L^{p-1}\psi \, dx\,dy$$ 
for every $\psi\in H^1(D(L))$. Choosing $v = v_L$ in the variational inequality and substituting $\psi = v_L-u_L$ into the nonlinear term, we obtain:
\begin{align*}
0 &\leq S_{\omega,L}'(u_L)[v_L-u_L] = a_{\omega,L}(u_L,v_L-u_L) - \int_{D(L)} u_L^{p-1}(v_L-u_L) \, dx\,dy \\
&= a_{\omega,L}(u_L,v_L-u_L) - a_{\omega,L}(v_L,v_L-u_L) = -a_{\omega,L}(u_L-v_L,u_L-v_L).
\end{align*}
Since $\omega>0$, the bilinear form $a_{\omega,L}$ is strictly coercive on $H^1(D(L))$, which forces $u_L = v_L$. Substituting $v_L=u_L$ back into the resolvent definition yields the unconstrained weak formulation.
\end{proof}

The next part of this section shows that \eqref{eqResolventInv} holds. The next lemma highlights the importance of the \textit{diagonal} symmetry $\tau$ inside the definition of $\mathcal{K}(L)$. 

\begin{lemma}[Regularity supplied by diagonal symmetry]
\label{lemDiagRegular}
Let \(F\in L^2(D(L))\) be invariant under the two coordinate
reflections and under the diagonal symmetry \(\tau\). If
\[
V=\mathcal R_{\omega,L}F,
\]
then
\[
V|_{D_+(L)}\in H^2(D_+(L)).
\]
Moreover,
\[
\|V\|_{H^2(D_+(L))}
\leq C_{\omega,L}\|F\|_{L^2(D(L))}.
\]
\end{lemma}

\begin{proof}
Let $g\in\{\sigma_x,\sigma_y,\tau\}$, where $\sigma_x(x,y)=(-x,y)$
and $\sigma_y(x,y)=(x,-y)$, and set $W\coloneqq V\circ g$. Since $g$
is an orthogonal transformation preserving $D(L)$, for every
$\psi\in H^1(D(L))$ we have
\[
\int_{D(L)}
\left(\nabla W\cdot\nabla\psi+\omega W\psi\right)
=
\int_{D(L)}
\left(\nabla V\cdot\nabla(\psi\circ g)
+\omega V(\psi\circ g)\right).
\]
Using the weak equation for $V$ and changing variables, this becomes
\[
\int_{D(L)}F(\psi\circ g)
=
\int_{D(L)}(F\circ g)\psi
=
\int_{D(L)}F\psi,
\]
because $F\circ g=F$. Thus $W$ and $V$ solve the same Neumann
resolvent problem. By uniqueness,
\[
V\circ g=V.
\]
Applying this to $\sigma_x$, $\sigma_y$, and $\tau$ proves that $V$
has all symmetries of $F$.

 The domain \(\Pi(L)\) defined in \eqref{eqOctant} is a fundamental domain for
these symmetries. Its boundary consists of
\[
\Gamma_0=\{(x,0):x>0\},\qquad
\Gamma_\Delta=\{(t,t):0<t<L/2\},\qquad
\Gamma_{L/2}=\{(x,L/2):x>L/2\}.
\]
For every \(\varphi\in H^1(\Pi(L))\), extend \(\varphi\) first to $D_+(L)$ by reflection and then to \(D(L)\)
by Lemma \ref{lemProjection}. Using the weak equation for \(V\) and
dividing by the number of copies gives
\[
\int_{\Pi(L)}
\left(\nabla V\cdot\nabla\varphi+\omega V\varphi\right)
=
\int_{\Pi(L)}F\varphi .
\]
Thus \(V|_{\Pi(L)}\) solves the Neumann problem on \(\Pi(L)\).

Now \(\Pi(L)\) is convex: its two finite angles are \(\pi/4\) and
\(3\pi/4\). Standard \(H^2\)-regularity for the Neumann problem on
convex polygonal domains (for instance \cite[Theorem 3.2.1.3]{Gr11}), together with the standard strip estimate
at infinity, yields
\[
V|_{\Pi(L)}\in H^2(\Pi(L)).
\]
Finally, the Neumann condition on \(\Gamma_\Delta\) permits even
reflection across the diagonal. This gives \(H^2\)-regularity on
\(D_+(L)\).
\end{proof}

\begin{remark}[The corner obstruction]
\label{rem:corner_obstruction}
Without the diagonal reflection assumption, Lemma \ref{lemDiagRegular} fails. Indeed, consider the corner in \(C=(L/2, L/2)\). Near \(C\), the first-quadrant
domain \(D_+(L)\) is a wedge of opening \(3\pi/2\). More precisely,
writing
\[
x-L/2=-r\sin\vartheta,\qquad
y-L/2=r\cos\vartheta,
\qquad
0<r<r_0,\quad 0<\vartheta<\frac{3\pi}{2},
\]
the two exterior boundary rays correspond to
\(\vartheta=0\) and \(\vartheta=3\pi/2\).

For the Neumann Laplacian on this wedge, the first angular
mode is
\[
\mathsf S_C(r,\vartheta)
\coloneqq
\chi(r)r^{2/3}\cos\left(\frac{2\vartheta}{3}\right),
\]
where \(\chi\) is a cut-off function equal to one near \(r=0\) \cite{Gr11}.
This function belongs to \(H^1\), but not to \(H^2\). Indeed, near
\(C\),
\begin{equation}\label{eqS_cDecay}
    |\nabla \mathsf S_C|\sim r^{-1/3},
\qquad
|D^2\mathsf S_C|\sim r^{-4/3},
\end{equation}
and hence
\[
\int_0^{r_0}r^{-2/3}r\,dr<\infty,
\qquad \mbox{ and } \quad
\int_0^{r_0}r^{-8/3}r\,dr=\infty.
\]
In particular, for \(\chi=1\),
\[
\partial_x\mathsf S_C
=
-\frac{2}{3}r^{-1/3}\sin\left(\frac{\vartheta}{3}\right),
\qquad
\partial_y\mathsf S_C
=
\frac{2}{3}r^{-1/3}\cos\left(\frac{\vartheta}{3}\right).
\]
Thus $
\partial_x\mathsf S_C,\ \partial_y\mathsf S_C
\in L^2 \setminus H^1.
$

Let \(F\in L^2(D(L))\), and let
\(V=\mathcal R_{\omega,L}F\). Standard Neumann corner regularity
theory gives the local decomposition
\[
V=V_{\rm reg}+c_C(F)\mathsf S_C
\qquad\text{near }C,
\]
where \(V_{\rm reg}\in H^2\) locally and \(c_C(F)\) is the leading
corner coefficient (see \cite[Section 8.2]{Gr11}). The two axial symmetries do not force
\(c_C(F)=0\), since they map \(C\) to other corners rather than fixing
it. Therefore, whenever \(c_C(F)\neq0\),
\[
V\notin H^2(D_+(L)\cap B_r(C)),\qquad
\partial_j V\notin H^1(D_+(L)\cap B_r(C))
\]
for small positive $r$. In particular, one gets from \eqref{eqS_cDecay} that
\[
\mathcal R_{\omega,L}:L^2(D(L))
\to H^{5/3-\varepsilon}(D(L))
\qquad\text{for every }\varepsilon>0,
\]
hence certainly not into \(H^2(D(L))\).
\end{remark}

The additional regularity from Lemma \ref{lemDiagRegular} lets us differentiate the resolvent identity.

\begin{lemma}[Differentiating resolvent identity]
\label{lem:D4_derivative_identity}
Let \(F\in H^1(D(L))\) be invariant under the two coordinate
reflections and under \(\tau\). Let
\(V=\mathcal R_{\omega,L}F\), and put \(W\coloneqq\partial_xV\).
Define
\[
\Gamma_D^x\coloneqq
\bigl(\{0\}\times(0,\infty)\bigr)
\cup
\bigl(\{L/2\}\times(L/2,\infty)\bigr),
\]
and
\[
\mathcal V_x\coloneqq
\left\{
\varphi\in H^1(D_+(L)):
\operatorname{Tr}\varphi=0
\text{ on }\Gamma_D^x
\right\}.
\]
Then \(W\in\mathcal V_x\), and, for every \(\varphi\in\mathcal V_x\),
\begin{equation}
    \label{eqDerEquat}
    \int_{D_+(L)}
\left(\nabla W\cdot\nabla\varphi+\omega W\varphi\right)
=
\int_{D_+(L)} (\partial_xF)\varphi
\end{equation}
\end{lemma}

\begin{proof}
By Lemma~\ref{lemDiagRegular}, the restricted resolvent solution satisfies \(V|_{D_+(L)}\in H^2(D_+(L))\), which ensures that \(W \coloneqq \partial_x V \in H^1(D_+(L))\). 

We first verify the boundary trace of \(W\). The domain \(D_+(L)\) has two vertical boundary components where the outward normal vector is purely horizontal. On the symmetry axis \(\{0\}\times(0,\infty)\), the evenness of \(V\) in the \(x\)-variable implies \(\partial_x V(0,y) = 0\). On the outer boundary \(\{L/2\}\times(L/2,\infty)\), the original Neumann boundary condition for \(V\) yields \(\partial_\nu V = \partial_x V = 0\). Together, these force \(W = 0\) on \(\Gamma_D^x\), proving that \(W \in \mathcal V_x\).

Let \(\varphi\) be a smooth test function vanishing on \(\Gamma_D^x\). Testing the weak equation for \(V\) on \(D_+(L)\) with \(\psi = -\partial_x\varphi\) yields:
\[
\int_{D_+(L)} \left( -\nabla V\cdot\partial_x(\nabla\varphi) - \omega V\partial_x\varphi \right) dx\,dy = -\int_{D_+(L)} F\,\partial_x\varphi \, dx\,dy.
\]
We integrate each term by parts with respect to \(x\). For the mass and source terms, the boundary integrals involve \(\varphi n_x\), where \(n_x\) is the \(x\)-component of the outward normal. Because \(n_x = 0\) on all horizontal boundaries and \(\varphi = 0\) on the vertical boundaries \(\Gamma_D^x\), these boundary terms vanish identically, yielding:
\[
-\int_{D_+(L)} \omega V\partial_x\varphi = \int_{D_+(L)} \omega (\partial_x V) \varphi = \int_{D_+(L)} \omega W \varphi, \quad \text{and} \quad -\int_{D_+(L)} F\partial_x\varphi = \int_{D_+(L)} (\partial_x F)\varphi.
\]
For the kinetic term, integrating by parts in \(x\) produces:
\[
-\int_{D_+(L)} \nabla V\cdot\partial_x(\nabla\varphi) = \int_{D_+(L)} \partial_x(\nabla V)\cdot\nabla\varphi - \int_{\partial D_+(L)} (\nabla V\cdot\nabla\varphi) n_x \, ds.
\]
We evaluate this boundary integral. On the horizontal boundary segments, \(n_x = 0\), so the integral is equal to zero. On the vertical boundary segments \(\Gamma_D^x\), the tangential derivative of the test function vanishes (\(\partial_y\varphi = 0\)) because \(\varphi = 0\) everywhere along these vertical lines. Furthermore, as established above, \(\partial_x V = W = 0\) on \(\Gamma_D^x\). Therefore, 
\[
(\nabla V\cdot\nabla\varphi) n_x = (\partial_x V)(\partial_x \varphi) + (\partial_y V)(\partial_y \varphi) = (0)(\partial_x \varphi) + (\partial_y V)(0) = 0.
\]
Thus, the boundary integral vanishes, leaving 
\[-\int_{D_+(L)} \nabla V\cdot\partial_x(\nabla\varphi) = \int_{D_+(L)} \nabla W \cdot\nabla\varphi.\]

Combining these terms yields \eqref{eqDerEquat} for smooth test functions. Since smooth functions vanishing on \(\Gamma_D^x\) are dense in \(\mathcal V_x\), the identity extends by continuity to every \(\varphi \in \mathcal V_x\).
\end{proof}

Since we can differentiate the resolvent, we can now show \eqref{eqResolventInv}.

\begin{lemma}[Resolvent invariance on the cone]
\label{lem:cone_invariance}
Property \eqref{eqResolventInv} holds; that is for every \(u\in\mathcal K(L)\),
\[
\mathcal R_{\omega,L}(u^{p-1})\in\mathcal K(L).
\label{eq:resolvent_invariance}
\]
\end{lemma}

\begin{proof}
For \(m>0\), put
\[
T_m(t)\coloneqq\min\{t,m\},\qquad
f_m\coloneqq T_m(u)^{p-1},\qquad
v_m\coloneqq\mathcal R_{\omega,L}f_m.
\]
Note that \(f_m\in H^1(D(L))\), it has the same symmetries as \(u\),
and
\[
\partial_jf_m
=
(p-1)T_m(u)^{p-2}\mathbf 1_{\{u<m\}}\partial_ju
\leq0
\quad\text{a.e.\ in }D_+(L),
\qquad j\in\{x,y\}.
\]
By uniqueness, \(v_m\) has the same symmetries. Set
\(w_m\coloneqq\partial_xv_m\). Using
Lemma~\ref{lem:D4_derivative_identity} with
\(\varphi=w_m^+\) gives
\begin{equation}
    \label{eqDerEq2}
    \int_{D_+(L)}
\left(
|\nabla w_m^+|^2+\omega|w_m^+|^2
\right) =
-\int_{D_+(L)}f_m\,\partial_xw_m^+ =
\int_{D_+(L)}(\partial_xf_m)w_m^+
\leq0.
\end{equation}
Hence \(\partial_xv_m\leq0\) a.e. in \(D_+(L)\). Interchanging \(x\)
and \(y\) gives \(\partial_yv_m\leq0\).

Moreover, testing the resolvent equation with
\(v_m^-=\max\{-v_m,0\}\) gives
\[
-\|v_m^-\|_{H^1_\omega(D(L))}^2
=
\int_{D(L)}f_m v_m^-\geq0,
\]
so \(v_m\geq0\). Therefore \(v_m\in\mathcal K(L)\).

Finally, \(u^{p-1}\in L^2(D(L))\), and dominated convergence gives
\[
f_m\to u^{p-1}
\qquad\text{in }L^2(D(L)).
\]
Writing \(v=\mathcal R_{\omega,L}(u^{p-1})\), coercivity yields
\[
\|v_m-v\|_{H^1_\omega(D(L))}
\leq
\omega^{-1/2}\|f_m-u^{p-1}\|_{L^2(D(L))}
\to0.
\]
Since \(\mathcal K(L)\) is strongly closed in \(H^1(D(L))\),
we obtain \(v\in\mathcal K(L)\).
\end{proof}

In the rest of this section, we discuss the roles of monotonicity and the two symmetry assumptions that characterize the cone $\mathcal{K}(L)$. First, we discuss the role of the diagonal symmetry. Without it, one obtains \eqref{eqDerEquat} only formally, since the pairing cannot be justified, as Remark \ref{rem:corner_obstruction} shows. And indeed the resolvent property \eqref{eqResolventInv} fails with only axis symmetry, as shown in the following proposition.
 
\begin{proposition}[Failure of resolvent invariance without diagonal symmetry]
\label{prop:axial_resolvent_failure}
Let 
\begin{equation*}
     \mathcal K_{\rm ax}(L)\coloneqq
\left\{
u\in H^1(D(L)):
\begin{array}{l}
u(-x,y)=u(x,y)=u(x,-y),\\
\partial_xu\leq0,\quad \partial_yu\leq0
\quad\text{a.e.\ in }D_+(L)
\end{array}
\right\}.
\end{equation*}
There exists $u\in\mathcal K_{\rm ax}(L)$ such that
\[
\mathcal R_{\omega,L}(u^{p-1})\notin\mathcal K_{\rm ax}(L).
\]
\end{proposition}

\begin{proof}
Let $C=(L/2,L/2)$.  Let $S_{\omega,C}$ be the local Neumann singular
solution of $(-\Delta+\omega)S_{\omega,C}=0$ near $C$, normalized by
\[
S_{\omega,C}(r,\vartheta)
=r^{2/3}\Phi(\vartheta)+O(r^{8/3}),
\]
where $\Phi(\vartheta) = \cos(2\vartheta/3)$ is the leading angular mode.  It is odd under the
diagonal reflection and, on a sufficiently small closed angular
sub-sector, satisfies, for some $\kappa>0$,
\[
\partial_xS_{\omega,C}\leq-\kappa r^{-1/3},
\qquad
\partial_yS_{\omega,C}\geq\kappa r^{-1/3}.
\]

Choose a radial cutoff $\chi$ supported in a small neighborhood of
$C$, with $\chi=1$ near $C$.  Extend
\[
\Psi\coloneqq \chi S_{\omega,C}
\]
evenly with respect to the two coordinate reflections. Then $\Psi$
satisfies the Neumann condition, is axially symmetric and diagonal
antisymmetric.  Since $S_{\omega,C}$ solves the homogeneous equation
near $C$,
\[
G\coloneqq(-\Delta+\omega)\Psi
\]
is smooth and supported in the cutoff annulus, away from the corners.
Consequently,
\[
\mathcal R_{\omega,L}G=\Psi,
\qquad
c_C(G)=1.
\]

Fix $a>0$ and define
\[
F_0(x,y)\coloneqq e^{-a(|x|+|y|)}.
\]
This function is axially and diagonally symmetric, positive, and
strictly decreasing in both variables on $D_+(L)$.  Hence
$c_C(F_0)=0$.  For sufficiently small $\varepsilon\neq0$, the function
\[
F_\varepsilon\coloneqq F_0+\varepsilon G
\]
remains positive, axially symmetric, and satisfies
\[
\partial_xF_\varepsilon\leq0,
\qquad
\partial_yF_\varepsilon\leq0
\quad\text{a.e.\ in }D_+(L).
\]
Indeed, $G$ and $\nabla G$ are bounded and compactly supported away
from the axes, whereas $F_0$ and $-\partial_jF_0$ have a strictly
positive lower bound on that compact support.

Set
\[
u_\varepsilon\coloneqq F_\varepsilon^{1/(p-1)}.
\]
Then $u_\varepsilon\in H^1(D(L))$, it has the axial symmetries, and
\[
\partial_j u_\varepsilon
=
\frac1{p-1}F_\varepsilon^{(2-p)/(p-1)}
\partial_jF_\varepsilon\leq0.
\]
Thus $u_\varepsilon\in\mathcal K_{\rm ax}(L)$ and
$u_\varepsilon^{p-1}=F_\varepsilon$. Moreover,
\[
c_C(F_\varepsilon)=\varepsilon.
\]
Let
\[
V_\varepsilon
=\mathcal R_{\omega,L}(u_\varepsilon^{p-1}).
\]
Suppose first that $\varepsilon>0$ and
$V_\varepsilon\in\mathcal K_{\rm ax}(L)$. Near $C$,
\[
V_\varepsilon
=V_{\rm reg}+\varepsilon\mathsf S_C,
\qquad
V_{\rm reg}\in H^2(\mathcal U_C\cap D_+(L)).
\]
Since $\partial_yV_\varepsilon\leq0$, on the chosen angular sub-sector
we obtain
\[
\partial_yV_{\rm reg}
\leq-\varepsilon\kappa r^{-1/3}.
\]
Hence
\[
\int |\partial_yV_{\rm reg}|^6
\geq C\int_0^\rho r^{-2}\,r\,dr=\infty.
\]
This is impossible because
\[
V_{\rm reg}\in H^2
\quad\to\quad
\partial_yV_{\rm reg}\in H^1\hookrightarrow L^6
\]
on the bounded corner sector.  If $\varepsilon<0$, the same argument
uses $\partial_x\mathsf S_C<0$.  Therefore
\[
\mathcal R_{\omega,L}(u_\varepsilon^{p-1})
\notin\mathcal K_{\rm ax}(L).
\]
\end{proof}

Next, we present an argument for why the diagonal symmetry $\tau$ alone is insufficient. Indeed, the diagonal reflection fixes only the two corners in the first and third quadrants, and bisects their re-entrant angles. Consequently, alone, this symmetry does not guarantee global
\(H^2(D(L))\) regularity, since the corner singularity described in
Remark~\ref{rem:corner_obstruction} can persist at the two off-diagonal corners.

Finally, without assuming monotonicity and symmetries, the minimizing problem admits escaping solitons, and compactness of the minimizing sequence is difficult to show for small values of $L$, as we show in this remark. 

\begin{remark}[Runaway sequences in the unrestricted problem]
\label{rem:runaway_sequences}
The unrestricted variational problem on $D(L)$ suffers from translation
loss of compactness along the four infinite arms. Indeed, let
$\phi_\omega$ be the one-dimensional soliton on $\mathbb R$, and let
$\eta\in C^\infty([0,\infty))$ satisfy $\eta=0$ near $0$ and
$\eta=1$ on $[1,\infty)$. On the right horizontal arm, with
$s=x-L/2$, define
\[
\psi_n(x,y)\coloneqq\eta(s)\phi_\omega(s-n),
\]
and extend $\psi_n$ by zero to the remainder of $D(L)$. Then
\[
\|\psi_n\|_{H^1_\omega(D(L))}^2
\to
L\|\phi_\omega\|_{H^1_\omega(\mathbb R)}^2,
\qquad
\|\psi_n\|_{L^p(D(L))}^p
\to
L\|\phi_\omega\|_{L^p(\mathbb R)}^p.
\]
After the scaling $t_n = \pi (\psi_n)>0$ where $\pi$ is defined in \eqref{eqNehariProj}, one has
$t_n\psi_n\in\mathcal N^{\rm free}(L)$, $t_n\to1$, and
\[
S_{\omega,L}(t_n\psi_n)\to
L s_{\omega,\mathbb R}.
\]
Moreover, $t_n\psi_n\rightharpoonup0$ in $H^1(D(L))$, while its
$L^p$ norm does not vanish.

The runaway mechanism is not merely a technical lack of compactness.
Indeed, Proposition \ref{prop:Identification} gives
\[
\frac{s_{\omega}(L)}{L}
\to
s_{\omega,\mathcal{G}}=2s_{\omega,\mathbb R},
\]
where
\[
s_{\omega,\mathbb R}
=
S_{\omega,\mathbb R}(\phi_\omega),
\]
where $\phi_\omega$ is the positive one-dimensional soliton defined
in \eqref{eqSechProf}. The escaping sequence above provides an upper bound for $s_{\omega}^{\rm free}(L)$ where
\[
\mathcal N^{\mathrm{free}}(L)
=\{u\in H^1(D(L);\mathbb R)\setminus\{0\}:I_{\omega,L}(u)=0\},
\qquad
s_\omega^{\mathrm{free}}(L)
=\inf_{u\in\mathcal N^{\mathrm{free}}(L)}S_{\omega,L}(u).
\]
to the minimal action level, as in particular 
\[
s_{\omega}^{\rm free}(L)
\leq Ls_{\omega,\mathbb R},
\]
Consequently, for all sufficiently small $L$,
\[
s_{\omega}^{\rm free}(L)
<
s_{\omega}(L).
\]
Thus, the symmetric solution constructed in this Section cannot be
an unrestricted ground state in the shrinking regime, as a profile
escaping along a single arm is energetically preferable. 
\end{remark}

\subsection{Variational characterization on the star graph}
\label{sec1DProblem}

In this subsection, we establish that the ground state profile $\phi_{\omega,\mathcal{G}}$ defined in \eqref{eqGraphGS} for $a_i = 0$ for any $i = 1,...,4$ can be characterized as the unique constrained minimizer of the action functional $S_{\omega, \mathcal{G}}$ on the Nehari manifold within a symmetric, monotonic subset of $H^1(\mathcal{G})$. We define the set of symmetric, monotonically decreasing functions on $\mathcal{G}$ by
\begin{equation}
    K_{sym}(\mathcal{G}) \coloneqq \left\{ v = (v_i)_{i=1,...,4} \in H^1(\mathcal{G}) : 
    \begin{array}{l}
        v_1(x) = v_2(x) = v_3(x) = v_4(x) \text{ for a.e. } x \in [0, \infty), \\
        v_1'(x) \le 0 \text{ for a.e. } x \in [0, \infty)
    \end{array}
    \right\}.
    \label{eq:K_graph}
\end{equation}
The Nehari manifold restricted to this set is given by
\begin{equation}
    \mathcal{N}(\mathcal{G}) \coloneqq \left\{ v \in K_{sym}(\mathcal{G}) \setminus \{0\} : I_{\omega, \mathcal{G}}(v) := \langle S'_{\omega, \mathcal{G}}(v), v \rangle =  0 \right\},
    \label{eq:nehari_graph}
\end{equation}
and the ground state energy level is defined by 
\begin{equation}
    \label{eq:min_graph}
    s_{\omega, \mathcal{G}} \coloneqq \inf_{v \in \mathcal{N}(\mathcal{G})} S_{\omega, \mathcal{G}}(v).
\end{equation}
\begin{lemma}[Compactness on the star graph]
\label{lem:compactness_graph}
For every $q > 2$, the embedding $K_{sym}(\mathcal{G}) \hookrightarrow L^q(\mathcal{G})$ is compact.
\end{lemma}

\begin{proof}
For a bounded sequence $(v_n) \subset K_{sym}(\mathcal{G})$, the edge-wise symmetry and monotonicity $v_{n,1} \ge 0$ and $v_{n,1}' \le 0$ imply that 
$$x |v_{n,1}(x)|^2 \le \int_0^x |v_{n,1}(y)|^2 dy \le \|v_{n,1}\|_{L^2}^2.$$
This yields the uniform pointwise tail estimate 
$$|v_{n,1}(x)|^2 \le \frac{1}{x} \|v_n\|_{H^1(\mathcal{G})}^2$$
for $x > 0$. This guarantees uniform decay at infinity, which, combined with the local Rellich compactness theorem, proves the statement.
\end{proof}

\begin{theorem}[Characterization of the graph minimizer]
\label{thm:existence_graph}
For every $\omega > 0$ and $p > 2$, the minimization problem \eqref{eq:min_graph} admits $\Phi_{\omega,\mathcal{G}}$ as its unique positive minimizer.
\end{theorem}

\begin{proof}
The existence of a non-trivial constrained minimizer $V \in \mathcal{N}(\mathcal{G}) \cap K_{sym}(\mathcal{G})$ follows from Lemma \ref{lem:compactness_graph}. To show that $V$ is a true critical point of the action $S_{\omega, \mathcal{G}}$ over the full space $H^1(\mathcal{G})$, we verify that the resolvent operator $\mathcal{R}_{\omega, \mathcal{G}} \coloneqq (-\Delta_{\mathcal{G}} + \omega)^{-1}$ preserves the cone $K_{sym}(\mathcal{G})$, following a similar, and actually easier, approach to that in Section \ref{sec:resolvent}.

Let $f \in K_{sym}(\mathcal{G})$ and define $W = (w_1, w_2, w_3, w_4) \coloneqq \mathcal{R}_{\omega, \mathcal{G}}(f^{p-1})$  as the unique weak solution to $$-\Delta_{\mathcal{G}} W + \omega W = f^{p-1}$$
in $H^1(\mathcal{G})$. By the symmetry of  $f$ it follows that $w_1 = w_2 = w_3 = w_4 \coloneqq w$.

Since $f \in H^1(\mathcal{G})$, $f^{p-1}$ belongs to $H^1(\mathcal{G})$ as $H^1$ is an algebra in one dimension. Standard one-dimensional elliptic regularity on each edge yields $W \in \bigoplus_{i=1}^4 H^3(e_i)$. This allows us to integrate by parts in the weak equation defining the resolvent:
\[
\int_{\mathcal{G}} (\nabla W \cdot \nabla \varphi + \omega W \varphi) dx = \sum_{i=1}^4 \int_{0}^\infty (-w_i'' + \omega w_i)\varphi_i dx - \sum_{i=1}^4 w_i'(0)\varphi_i(0) = \int_{\mathcal{G}} f^{p-1}\varphi dx.
\]
Because the strong equation $$-w_i'' + \omega w_i = f_i^{p-1}$$
holds almost everywhere on each edge, the volume integrals cancel, leaving the boundary term 
$$\varphi(0) \sum_{i=1}^4 w_i'(0) = 0$$
for all test functions $\varphi \in H^1(\mathcal{G})$. This establishes the Kirchhoff condition 
$$\sum_{i=1}^4 w_i'(0) = 0.$$ 
The symmetry forces $4w'(0) = 0$, so $w'(0) = 0$.

Differentiating the equation for $w$ on $(0, \infty)$ yields 
\[
-(w')'' + \omega (w') = (p-1)f_1^{p-2}f_1' \le 0 \quad \text{on } (0, \infty).
\]
We test this equation against $(w')^+ \coloneqq \max\{w', 0\}$. Because $w'(0) = 0$, the test function $(w')^+$ has zero trace at the origin, so the boundary evaluation upon integration by parts vanishes. We obtain:
\[
\int_0^\infty
\left|\partial_x(w')^+\right|^2+\omega\left|(w')^+\right|^2 dx
=\int_0^\infty (p-1)f_1^{p-2}f'_1(w')^+ dx \le0.
\]
The right-hand side is non-positive because $f_1' \le 0$ by the assumption that $f \in K_{sym}(\mathcal{G})$. Thus, $(w')^+ = 0$, meaning $w' \le 0$. This confirms $W \in K_{sym}(\mathcal{G})$.

The invariance of $K_{sym}(\mathcal{G})$ under the resolvent (as condition \eqref{eqResolventInv} for the cross) ensures that the constrained minimizer is a global weak solution to the nonlinear Schrödinger equation on $\mathcal{G}$. The invariance of $K_{sym}(\mathcal{G})$ under the resolvent (analogous to condition \eqref{eqResolventInv} for the cross) ensures that the constrained minimizer $V$ is a global weak solution to the nonlinear Schrödinger equation on $\mathcal{G}$, following the exact same coercivity argument as in Proposition \ref{prop:natural_cone}.

The limiting profile is nonnegative and nonzero. By the maximum principle, it is strictly positive in the interior of any half-line. The unique positive $H^1(0, \infty)$ solution satisfying $v'(0) = 0$ is the unshifted restriction of the even sech soliton, $\phi_\omega|_{[0,\infty)}$. The vertex continuity condition \eqref{eqGraphConti} is satisfied since all components evaluate to $\phi_\omega(0)$. Consequently, $V \equiv \Phi_{\omega,\mathcal{G}}$.
\end{proof}

\section{Properties of minimizers}
\label{sec:properties}

Having established the existence of a minimizer $u_L \in \mathcal{N}(L)$, we now show that it decays exponentially along the arms at a rate independent of the strip width $L$. In this section, we fix $\omega >0$
and we recall that for every $z\in H^1(D(L))$

\begin{equation*}
    \|z\|_{H^1_\omega(D(L))}^2  = \int_{D(L)} |\nabla z|^2 + \omega |z|^2 dz.
\end{equation*}
We start with the following scaling lemma.

\begin{lemma}[Scale-sensitive Sobolev estimates]
\label{lem:thin_Sobolev}
For every \(p>2\), there exists \(C=C(p,\omega)>0\), independent of
\(L\in(0,1]\), such that
\begin{equation}
\label{eq:thin_Sobolev}
\|z\|_{L^p(D(L))}^p
\leq
C L^{1-p/2}\|z\|_{H^1_\omega(D(L))}^p
\qquad\text{for every }z\in H^1(D(L)).
\end{equation}
\end{lemma}

\begin{proof}
Consider first the horizontal strip $h(L)$ defined in \eqref{eqHandV} and, for \(z\in H^1(h(L))\), set
\[
\widetilde z(s,\eta)\coloneqq z(s,L\eta),
\qquad
(s,\eta)\in\Sigma\coloneqq\mathbb R\times(-1/2,1/2).
\]
The Sobolev embedding on the fixed strip \(\Sigma\) gives
\[
\|\widetilde z\|_{L^p(\Sigma)}^p
\leq C_p\|\widetilde z\|_{H^1(\Sigma)}^p.
\]
Moreover,
\begin{align*}
\|\widetilde z\|_{H^1(\Sigma)}^2
&=
L^{-1}\int_{h(L)}
\left(|\partial_xz|^2+|z|^2\right)\,dx\,dy
+
L\int_{h(L)}|\partial_yz|^2\,dx\,dy \leq
C_\omega L^{-1}\|z\|_{H^1_\omega(D(L))}^2.
\end{align*}
Therefore,
\[
\|z\|_{L^p(h(L))}^p
=
L\|\widetilde z\|_{L^p(\Sigma)}^p
\leq
C L^{1-p/2}\|z\|_{H^1_\omega(D(L))}^p.
\]
The same estimate holds on the vertical strip \( v(L)\). Since \( D(L)=h(L)\cup v(L), \)
the conclusion follows.
\end{proof}

\begin{lemma}[Energy scale and nonvanishing]
\label{lem:energy_nonvanishing}
Let \(u_L\) be a minimizer to \eqref{eqMinProb}. Then there exist constants \(0<c\leq C\), independent of \(L\in(0,1]\),
such that
\begin{equation}
\label{eqTwoSideEnergy}
cL\leq  \|u_L\|_{H^1_\omega(D(L))}^2\leq CL.
\end{equation}
Moreover,
\begin{equation}
\label{eqLimsupEneGraph}
\limsup_{L\to0}\frac{s_{\omega}(L)}{L}
\leq s_{\omega,\mathcal G}.
\end{equation}
\end{lemma}

\begin{proof}
Let
\[
\Phi_{\omega,L}
\coloneqq
\mathcal C_L\bigl(\phi_\omega|_{\mathbb R^+}\bigr)
\]
be the cross-extension of the line ground state given by \eqref{eqCrossExte}. 
This function belongs to \(\mathcal K(L)\). Since
\(\phi_\omega'(0)=0\) and
\(\phi_\omega(0)^{p-2}=p\omega/2\), one has
\[
I_{\omega,L}(\Phi_{\omega,L})
=
L^2\bigl(\omega\phi_\omega(0)^2-\phi_\omega(0)^p\bigr)<0.
\]
The relaxed characterization of \(s_{\omega}(L)\) therefore gives
\begin{align*}
s_{\omega}(L)
&\leq
\frac{p-2}{2p}\|\Phi_{\omega,L}\|_{L^p(D(L))}^p =
L s_{\omega,\mathcal G}
+
\frac{p-2}{2p}L^2\phi_\omega(0)^p.
\end{align*}
This proves \eqref{eqLimsupEneGraph} and, by the Nehari identity,
the upper bound in \eqref{eqTwoSideEnergy}.

For the lower bound, the Nehari identity and
Lemma~\ref{lem:thin_Sobolev} yield
\[
 \|u_L\|_{H^1_\omega(D(L))}^2
=
\|u_L\|_{L^p(D(L))}^p
\leq
C L^{1-p/2} (\|u_L\|_{H^1_\omega(D(L))}^2)^{p/2}.
\]
Since \( \|u_L\|_{H^1_\omega(D(L))}^2>0\), division by \( \|u_L\|_{H^1_\omega(D(L))}^2\) gives
\[
 (\|u_L\|_{H^1_\omega(D(L))}^2)^{(p-2)/2}
\geq
C^{-1}L^{(p-2)/2},
\]
which proves the lower bound in \eqref{eqTwoSideEnergy}.
\end{proof}

\begin{proposition}[Uniform exponential decay]\label{prop:decay}
There exist constants $M_0 > 0$ and $R = R(\omega,p) > 0$, independent of $L$, such that for any $L \in (0, 1]$ and $(x,y) \in D(L)$ with $|(x,y)| \geq R$, it holds
$$
u_L(x,y) \le M_0 e^{-\sqrt{\omega/2} |(x,y)|}.
$$
\end{proposition}

\begin{proof}
By spatial symmetry, it suffices to prove the decay on the right horizontal arm $h^+(L)$ defined in \eqref{eqH+}. Using the pointwise tail estimate \eqref{eqPolyDecay} and substituting the energy bound \eqref{eqTwoSideEnergy}, we obtain for $x > L/2$:
\[
u_L(x,y)^2 \leq C_\omega \frac{\|u_L\|_{H^1_\omega(D(L))}^2}{x - L/2} \left(\frac{1}{L} + 2\right) \lesssim \frac{1 + 2L}{x - L/2}.
\]
Thus, $u_L$ decays to zero uniformly in $L$. We choose a sufficiently large radius $R \ge 1$, independent of $L$, so that for all $x \ge R$,
$$u_L(x,y)^{p-2} \le \frac{\omega}{2}.$$ 
On $(R, \infty) \times (-L/2, L/2)$, the positive solution $u_L$ satisfies 
$$-\Delta u_L + \frac{\omega}{2} u_L = u_L \left( u_L^{p-2} - \frac{\omega}{2} \right) \le 0.$$
We introduce the comparison function $w(x,y) = M e^{-\delta x}$ with $\delta = \sqrt{\omega / 2}$. By choosing 
$$M = e^{\delta R} \sqrt{C_1 / (R - 1/2)} \coloneqq M_0,$$
we guarantee that $w(R, y) \ge u_L(R, y)$. The maximum principle applied to $-\Delta + \omega/2$ forces $u_L \le w$ everywhere in $(R, \infty) \times (-L/2, L/2)$. 
\end{proof}

\begin{lemma}[Uniform Hessian estimate]
\label{lem:uniform_D4_H2}
Let $\omega>0$, $L>0$, and let $F\in L^2(D(L))$ be invariant
under the two coordinate reflections and the diagonal reflection
$\tau$. If $U=\mathcal R_{\omega,L}F$, then $U\in H^2(D(L))$ and
\[
\|D^2U\|_{L^2(D(L))}
\le \|\Delta U\|_{L^2(D(L))}
\le \|F\|_{L^2(D(L))}.
\]
Here $|D^2U|^2=\sum_{i,j=1}^2|\partial_{ij}U|^2$.
In particular,
\[
\|U\|_{H^2(D(L))}\le C_\omega\|F\|_{L^2(D(L))},
\]
with $C_\omega$ independent of $L \in (0,1]$.
\end{lemma}

\begin{proof}
By uniqueness, $U$ inherits the symmetries of $F$.
Lemma~\ref{lemDiagRegular}, followed by even reflection across
the coordinate axes, gives $U\in H^2(D(L))$.
Its restriction to  \(\Pi(L)\) defined in \eqref{eqOctant} satisfies homogeneous Neumann conditions.

We use the homogeneous Neumann estimate on bounded convex domains:
\begin{equation}
\label{eq:convex_neumann_hessian}
\|D^2z\|_{L^2(\Omega)}
\le \|\Delta z\|_{L^2(\Omega)},
\qquad
z\in H^2(\Omega),\quad \partial_\nu z=0.
\end{equation}
For instance, this follows from
\cite[Theorem~2.2, equation~(2.8) with $B(t)=t^2/2$ and $H(\xi)=|\xi|$]{ACCFM25}.
To apply this estimate on the unbounded chamber, choose
$\chi\in C^\infty([0,\infty))$ with $0\le\chi\le1$,
$\chi=1$ on $[0,1]$, and $\chi=0$ on $[2,\infty)$.
For $R>L$, set
\[
\Omega_R=\Pi(L)\cap\{x<3R\},
\qquad U_R(x,y)=\chi(x/R)U(x,y).
\]
The domain $\Omega_R$ is bounded and convex, and $U_R$ satisfies
homogeneous Neumann conditions on its entire boundary.
Indeed, the cutoff is constant near the diagonal boundary,
its normal derivative vanishes on the horizontal sides,
and $U_R$ vanishes near the artificial boundary $x=3R$.
Consequently, \eqref{eq:convex_neumann_hessian} yields
\[
\|D^2U_R\|_{L^2(\Omega_R)}
\le\|\Delta U_R\|_{L^2(\Omega_R)}.
\]
Regarded as a function on $\Pi(L)$ by zero extension past
$x=3R$, $U_R$ converges to $U$ in $H^2(\Pi(L))$ as
$R\to\infty$. This follows from $U\in H^2(\Pi(L))$
and by the fact that the cutoff derivative terms are bounded by
\[
\frac{C}{R}\|\nabla U\|_{L^2(\Pi(L))}
+\frac{C}{R^2}\|U\|_{L^2(\Pi(L))}.
\]
Passing to the limit gives
\[
\|D^2U\|_{L^2(\Pi(L))}
\le\|\Delta U\|_{L^2(\Pi(L))}.
\]
The cross consists of eight reflected copies of $\Pi(L)$.
Since orthogonal transformations preserve the Frobenius norm
of the Hessian and the Laplacian, symmetry gives
\[
\|D^2U\|_{L^2(D(L))}
\le\|\Delta U\|_{L^2(D(L))}.
\]
Finally, using $F=-\Delta U+\omega U$  we obtain
\[
\|F\|_2^2
=\|\Delta U\|_2^2
 +2\omega\|\nabla U\|_2^2
 +\omega^2\|U\|_2^2.
\]
This proves $\|\Delta U\|_2\le\|F\|_2$ and, together with
the Hessian estimate, the asserted full $H^2$ bound.
\end{proof}

\begin{proposition}[Uniform \(H^2\) and \(L^\infty\) bounds]
\label{prop:qualitative_cross}
Let \(u_L\in\mathcal K(L)\) be the  minimizer given by Theorem~\ref{thmExistMinim}. Then
\[
u_L\in H^2(D(L)).
\]
Moreover, there exists \(C>0\), independent of \(L\in(0,1]\), such that
\begin{equation}
\label{eq:uniform_H2}
\|u_L\|_{H^2(D(L))}\leq C L^{1/2},
\end{equation}
and
\begin{equation}
\label{eqUnifLInfty}
\|u_L\|_{L^\infty(D(L))}\leq C.
\end{equation}
\end{proposition}

\begin{proof}
Set \(F_L\coloneqq u_L^{p-1}\). Since \(u_L\in\mathcal K(L)\),
the function \(F_L\) has the same symmetries. By the
thin-domain Sobolev estimate with exponent \(2p-2\),
\[
\|u_L\|_{L^{2p-2}(D(L))}^{2p-2}
\leq
C L^{2-p}\|u_L\|_{H^1_\omega(D(L))}^{2p-2}.
\]
Using \eqref{eqTwoSideEnergy}, we obtain
\[
\|F_L\|_{L^2(D(L))}^2
=
\|u_L\|_{L^{2p-2}(D(L))}^{2p-2}
\leq
C L^{2-p}\bigl(CL\bigr)^{p-1}
\leq CL.
\]
Lemma~\ref{lem:uniform_D4_H2} therefore yields
\[
\|D^2u_L\|_{L^2(D(L))}
\leq C\|F_L\|_{L^2(D(L))}
\leq CL^{1/2}.
\]
Together with \eqref{eqTwoSideEnergy}, this proves
\eqref{eq:uniform_H2}.

We still need to prove the uniform \(L^\infty\)-bound. Consider first the strip $h(L)$ defined in \eqref{eqHandV}.
For \(z\in H^2(h(L))\), define
\[
\widetilde z(s,\eta)\coloneqq z(s,L\eta),
\qquad
(s,\eta)\in\mathbb R\times(-1/2,1/2).
\]
The Sobolev embedding on the fixed-width strip gives
\[
\|z\|_{L^\infty(h(L))}
=
\|\widetilde z\|_{L^\infty(h(1))}
\leq
C\|\widetilde z\|_{H^2(h(1))}
\leq
CL^{-1/2}\|z\|_{H^2(h(L))}.
\]
With the same estimate on the vertical strip $V(L)$, we obtain
\[
\|u_L\|_{L^\infty(D(L))}
\leq
CL^{-1/2}\|u_L\|_{H^2(D(L))}
\leq C,
\]
by \eqref{eq:uniform_H2}.
\end{proof}

We now show Theorem \ref{thm:existence}.

\begin{proof}[Proof of Theorem \ref{thm:existence}]
The existence follows from combining Theorem \ref{thmExistMinim} with Lemma \ref{lem:cone_invariance} and Proposition \ref{prop:natural_cone}. Positivity is proven in Proposition \ref{thmExistMinim} and Remark \ref{rem:StricPos} while the other properties are shown in Section \ref{sec:properties}
    
\end{proof}

\section{The Shrinking Limit}
\label{sec:shrinking_limit}

Let us define the scaled (independent of $L$) domains:
\[
I\coloneqq(-1/2,1/2),
\qquad
\Sigma\coloneqq\mathbb R^+\times I,
\qquad
K\coloneqq[-1/2,1/2]^2.
\]
For the right horizontal arm, we define
\[
w_L(s,\eta)
\coloneqq
u_L(s+L/2,L\eta),
\qquad
(s,\eta)\in\Sigma,
\]
The rescaled functions on all four arms agree with \(w_L\), up to outward-coordinate conventions, by axial and diagonal symmetry. We define the average and the deviation from the average 
\[
a_L(s)\coloneqq\int_I w_L(s,\eta)\,d\eta,
\qquad
\rho_L(s,\eta)\coloneqq w_L(s,\eta)-a_L(s).
\]
Note that we have 
\[
a_L\geq0,
\qquad
a_L'\leq0\quad\text{a.e.\ on }\mathbb R^+.
\]
Here, the monotonicity of \(a_L\) follows by averaging
\(\partial_xu_L\leq0\) on the right horizontal arm. Furthermore, in the central square, set
\[
v_L(\xi,\eta)\coloneqq u_L(L\xi,L\eta),
\qquad
(\xi,\eta)\in K.
\]

\begin{lemma}\label{lem:SomeConve}
The following hold:

\begin{equation}
\label{eq:average_energy_decomposition}
\begin{aligned}
\frac{ \|u_L\|_{H^1_\omega(D(L))}^2}{L}
={}&4\|a_L\|_{H^1_\omega(\mathbb R^+)}^2
+4\|\partial_s\rho_L\|_{L^2(\Sigma)}^2
+4L^{-2}\|\partial_\eta\rho_L\|_{L^2(\Sigma)}^2 \\
&+4\omega\|\rho_L\|_{L^2(\Sigma)}^2
+\frac1L\|\nabla v_L\|_{L^2(K)}^2
+\omega L\|v_L\|_{L^2(K)}^2 .
\end{aligned}
\end{equation}
In particular, the average is bounded while the deviation vanishes
\begin{equation}
\label{eqAverageUnifH1}
\|a_L\|_{H^1_\omega(\mathbb R^+)}\leq C,
\qquad
\|\rho_L\|_{L^2(\Sigma)}\leq CL.
\end{equation}
Finally, this scale-adjusted convergence holds
\begin{equation}
\label{eqAveNehariLim}
\left|
\frac{ \|u_L\|_{H^1_\omega(D(L))}^2}{L}
-
4\|a_L\|_{L^p(\mathbb R^+)}^p
\right|
\leq CL.
\end{equation}
\end{lemma}

\begin{proof}

The change of variables $x=s+L/2, y=L\eta $ on each arm, together with $ x=L\xi,  y=L\eta $ on the central square, gives
\begin{equation*}
\begin{aligned}
\frac{ \|u_L\|_{H^1_\omega(D(L))}^2}{L}
={}&4\int_\Sigma
\left(
|\partial_sw_L|^2
+
L^{-2}|\partial_\eta w_L|^2
+
\omega|w_L|^2
\right)\,ds\,d\eta \\
&+\frac1L\int_K|\nabla v_L|^2\,d\xi\,d\eta
+\omega L\int_K|v_L|^2\,d\xi\,d\eta.
\end{aligned}
\end{equation*}
Since
\begin{equation}
    \label{eqPoincareHyp}
    \int_I\rho_L(s,\eta)\,d\eta=0,
\end{equation}
one has
\[
\int_\Sigma|\partial_sw_L|^2
=
\int_0^\infty|a_L'|^2\,ds
+
\int_\Sigma|\partial_s\rho_L|^2\,ds\,d\eta,
\]
and likewise
\[
\int_\Sigma|w_L|^2
=
\int_0^\infty|a_L|^2\,ds
+
\int_\Sigma|\rho_L|^2\,ds\,d\eta.
\]
Moreover, \(\partial_\eta w_L=\partial_\eta\rho_L\). This proves
\eqref{eq:average_energy_decomposition}.

The bound for \(a_L\) follows from
\eqref{eqTwoSideEnergy} and
\eqref{eq:average_energy_decomposition}. Since \eqref{eqPoincareHyp}, the Poincar\'e--Wirtinger
inequality on \(I\) gives
\[
\|\rho_L\|_{L^2(\Sigma)}^2
\leq
C\|\partial_\eta\rho_L\|_{L^2(\Sigma)}^2
\leq CL^2,
\]
which proves \eqref{eqAverageUnifH1}.

By Proposition~\ref{prop:qualitative_cross} we get
$$\|w_L\|_{L^\infty(\Sigma)}+\|a_L\|_{L^\infty(\mathbb R^+)}\leq C.$$
Hence, it follows that
\[
\left|
\|w_L\|_{L^p(\Sigma)}^p
-
\|a_L\|_{L^p(\mathbb R^+)}^p
\right|
\leq
C\bigl(\|w_L\|_{L^2(\Sigma)}+\|a_L\|_{L^2(\mathbb R^+)}\bigr)
\|\rho_L\|_{L^2(\Sigma)}
\leq CL.
\]
The Nehari identity and the corresponding nonlinear scaling identity
\[
\frac{ \|u_L\|_{H^1_\omega(D(L))}^2}{L}
=
4\|w_L\|_{L^p(\Sigma)}^p
+
L\|v_L\|_{L^p(K)}^p
\]
then imply \eqref{eqAveNehariLim}, because
\(\|v_L\|_{L^p(K)}\leq C\) uniformly in $L$.
\end{proof}

\begin{lemma}[Compactness and nonvanishing of the arm profile]
\label{lem:common_arm_compactness}
For every sequence \(L_n\to0^+\), there are a subsequence, still denoted by $L_n$,
and a function \(\phi\in H^1(\mathbb R^+)\) such that
\[
a_{L_n}\rightharpoonup\phi
\quad\text{weakly in }H^1(\mathbb R^+),
\qquad
a_{L_n}\to\phi
\quad\text{strongly in }L^p(\mathbb R^+).
\]
Moreover,
\[
\phi\geq0,\qquad
\phi'\leq0\quad\text{a.e.},\qquad
\phi\not\equiv0.
\]
\end{lemma}

\begin{proof}
The \(H^1\)-bound follows from
\eqref{eqAverageUnifH1}. Thus, for any $(u_{L_n})$, we can take a weakly converging subsequence. Positivity and monotonicity are weakly
closed properties. For \(s>0\), monotonicity gives
\[
s\,a_L(s)^2
\leq
\int_0^s a_L(t)^2\,dt
\leq C.
\]
Consequently,
\[
\int_R^\infty a_L(s)^p\,ds
\leq
\sup_{s\geq R}a_L(s)^{p-2}
\int_R^\infty a_L(s)^2\,ds
\leq
CR^{-(p-2)/2}.
\]
Rellich compactness on \((0,R)\), followed by \(R\to\infty\), proves
strong \(L^p(\mathbb R^+)\)-convergence, as in the proof of Proposition \ref{prop:compactness_HL}.

 Now, by \eqref{eqTwoSideEnergy} and
\eqref{eqAveNehariLim}, we get
\[ 4\Vert{}a_L\Vert{}_{L^p(\R^+)}^p \ge \frac{\Vert{}u_L\Vert{}_{H^1_\omega}^2}{L} - CL \ge c - CL >0
\]
for all sufficiently small \(L\). Thus \(\phi\not\equiv0\).
\end{proof}

\begin{proposition}[Identification of the graph limit]
\label{prop:Identification}
One has
\begin{equation}
\label{eqEnergyLimit}
\lim_{L\to0}\frac{s_{\omega}(L)}{L}
=
s_{\omega,\mathcal G}.
\end{equation}
Moreover, if \(a_{L_n}\rightharpoonup\phi\) as in 
Lemma~\ref{lem:common_arm_compactness}, then
\[
\phi=\phi_\omega|_{\mathbb R^+}.
\]
where $\phi_\omega$ is the one dimensional ground state defined in \eqref{eqSechProf}.
\end{proposition}

\begin{proof}
Let
\[
\Phi\coloneqq(\phi,\phi,\phi,\phi)\in K_{\rm sym}(\mathcal G).
\]
By \eqref{eqAveNehariLim} and strong \(L^p\)-convergence, it follows that
\begin{equation*}
\lim_{n\to\infty}\frac{ \|u_{L_n}\|_{H^1_\omega(D(L))}^2}{L_n}
= 4\|\phi\|_{L^p(\mathbb R^+)}^p.
\end{equation*}
On the other hand, the nonnegative terms in
\eqref{eq:average_energy_decomposition} and weak lower
semicontinuity yield
\[
4\|\phi\|_{H^1_\omega(\mathbb R^+)}^2
\leq
\liminf_{n\to\infty}\frac{ \|u_{L_n}\|_{H^1_\omega(D(L))}^2}{L_n}
=
4\|\phi\|_{L^p(\mathbb R^+)}^p.
\]
Thus we obtain
\[
I_{\omega,\mathcal G}(\Phi)\leq0.
\]
Set
\[
t\coloneqq
\left(
\frac{\|\phi\|_{H^1_\omega(\mathbb R^+)}^2}
     {\|\phi\|_{L^p(\mathbb R^+)}^p}
\right)^{1/(p-2)}
\in(0,1].
\]
Then \(t\Phi\in\mathcal N(\mathcal G)\cap K_{\rm sym}(\mathcal G)\),
and hence
\begin{align*}
s_{\omega,\mathcal G} \leq
S_{\omega,\mathcal G}(t\Phi) =
\frac{p-2}{2p}\,
t^p\|\Phi\|_{L^p(\mathcal G)}^p \leq
\frac{p-2}{2p}\,
\|\Phi\|_{L^p(\mathcal G)}^p =
\lim_{n\to\infty}
\frac{s_{\omega}(L_n)}{L_n}.
\end{align*}
Together with \eqref{eqLimsupEneGraph}, this proves
\eqref{eqEnergyLimit}.

All inequalities in the preceding display must be equalities. Therefore
\(t=1\), \(\Phi\in\mathcal N(\mathcal G)\), and
\[
S_{\omega,\mathcal G}(\Phi)=s_{\omega,\mathcal G}.
\]
Since $\phi \geq 0$, the uniqueness property on the one-dimensional ground state stated in Theorem~\ref{thm:existence_graph} yields
\[
\Phi=\phi_{\omega,\mathcal G},
\]
which proves the asserted identification of \(\phi\).
\end{proof}

\begin{proposition}[Strong collapse on the fixed reference domains]
\label{prop:strong_D4_collapse}
As \(L\to0\),
\[
a_L\to\phi_\omega|_{\mathbb R^+}
\quad\text{strongly in }H^1(\mathbb R^+), \quad  w_L\to\phi_\omega(s)
\quad\text{strongly in }H^1(\Sigma),
\]
and
\[
v_L\to\phi_\omega(0)
\quad\text{strongly in }H^1(K).
\]
\end{proposition}

\begin{proof}
By Proposition~\ref{prop:Identification}, every sequence $L_n \to 0^+$ admits a subsequence, still denoted by $L_n$ such that $a_{L_n} \rightharpoonup \phi_\omega|_{\mathbb R^+}$ weakly. Hence the full
family converges weakly to this function in \(H^1(\mathbb R^+)\).

Using \eqref{eqEnergyLimit},
\eqref{eq:average_energy_decomposition}, and weak lower
semicontinuity, we obtain
$$
\|a_L\|_{H^1_\omega(\mathbb R^+)}
\to
\|\phi_\omega\|_{H^1_\omega(\mathbb R^+)},$$
which implies that \(a_L\to\phi_\omega\) strongly in \(H^1(\mathbb R^+)\).
Moreover, every nonnegative remainder term in
\eqref{eq:average_energy_decomposition} tends to zero. In particular,
\[
\|\partial_s\rho_L\|_{L^2(\Sigma)}\to0,
\qquad
L^{-1}\|\partial_\eta\rho_L\|_{L^2(\Sigma)}\to0,
\qquad
L^{-1/2}\|\nabla v_L\|_{L^2(K)}\to0.
\]
Together with \(\|\rho_L\|_{L^2(\Sigma)}\leq CL\), this proves
\[
\| w_L - a_L \|_{H^1(\Sigma)} \to 0
\]
and hence the asserted convergence of \(w_L\).

The remaining task is to identify the central limit. Let
\[
m_L\coloneqq\int_Kv_L(\xi,\eta)\,d\xi\,d\eta.
\]
The trace theorem and the Poincar\'e inequality on \(K\) imply
\[
|a_L(0)-m_L|
\leq
C\|\nabla v_L\|_{L^2(K)}.
\]
Since \(a_L(0)\to\phi_\omega(0)\), we obtain
\(m_L\to\phi_\omega(0)\). A further use of the Poincar\'e inequality
on \(K\) now gives
\[
\|v_L-\phi_\omega(0)\|_{H^1(K)}\to0.
\]
\end{proof}

To complete the proof of Theorem \ref{thm:main_shrinking_limit}, we show that, using the resolvent estimates, we obtain $H^2$ convergence.

\begin{proof}[Proof of Theorem \ref{thm:main_shrinking_limit}.]
We write $\phi=\phi_\omega|_{\mathbb R^+}$ and
$e_L=u_L-\Phi_{\omega,L}$.
We first obtain a normalized $H^1$ convergence.
By the same orthogonality argument as in
\eqref{eq:average_energy_decomposition} we get
\begin{align*}
\frac{\|e_L\|_{H^1_\omega(D(L))}^2}{L}
={}&4\|a_L-\phi\|_{H^1_\omega(\mathbb R^+)}^2
   +4\|\partial_s\rho_L\|_{L^2(\Sigma)}^2
   +4L^{-2}\|\partial_\eta\rho_L\|_{L^2(\Sigma)}^2\\
&+4\omega\|\rho_L\|_{L^2(\Sigma)}^2
   +L^{-1}\|\nabla v_L\|_{L^2(K)}^2
   +\omega L\|v_L-\phi(0)\|_{L^2(K)}^2.
\end{align*}
Every term on the right tends to zero by Lemma \ref{lem:SomeConve} and 
Proposition~\ref{prop:strong_D4_collapse}. 
Since $\omega>0$ is fixed, this gives
\begin{equation}
\label{eqNormH1Conv}
L^{-1/2}\|e_L\|_{H^1(D(L))}\to 0.
\end{equation}

We next use the equation to improve this convergence.
The comparison profile $\Phi_{\omega,L}$ belongs to $H^2(D(L))$ (as the first derivatives match because $\phi_{\omega}'(0)=0$) which eliminates the interface measures) and satisfies the homogeneous Neumann boundary condition.
Consequently, in the weak sense on $D(L)$,
\[
(-\Delta+\omega)\Phi_{\omega,L}
=\Phi_{\omega,L}^{p-1}+\kappa\mathbf{1}_{q(L)},
\qquad
\kappa=\omega\phi(0)-\phi(0)^{p-1}.
\]
Indeed, the one-dimensional equation holds on every arm,
whereas $\Phi_{\omega,L}=\phi(0)$ on $q(L)$.
There are no interface terms because the first derivatives match. Subtracting this identity from the equation for $u_L$ \eqref{eqNlsCross} gives
\[
e_L=\mathcal R_{\omega,L}F_L,
\qquad
F_L=u_L^{p-1}-\Phi_{\omega,L}^{p-1}
       -\kappa\mathbf{1}_{q(L)}.
\]

By Proposition~\ref{prop:qualitative_cross} and the definition
of $\Phi_{\omega,L}$, both functions are uniformly bounded
in $L^\infty(D(L))$.
Since $p>2$, the map $t\mapsto t^{p-1}$ is Lipschitz on the
common bounded interval containing their nonnegative values.
As $|q(L)|=L^2$, we obtain
\[
\|F_L\|_{L^2(D(L))}
\le C\|e_L\|_{L^2(D(L))}+|\kappa|L
=o(L^{1/2}).
\]
The source $F_L$ has the coordinate and diagonal symmetries.
Lemma~\ref{lem:uniform_D4_H2} therefore yields
\[
\|D^2e_L\|_{L^2(D(L))}
\le C\|F_L\|_{L^2(D(L))}=o(L^{1/2}).
\]
Together with \eqref{eqNormH1Conv}, this yields
\[
L^{-1/2}\|e_L\|_{H^2(D(L))}\to 0,
\]
which proves the theorem.
\end{proof}

\section*{Acknowledgments}

The author thanks Abbas Moameni for sharing with him his personal notes, which provided several useful points for the existence argument developed in this paper. The author thanks Stefan Le Coz for useful discussions about this work. The work of B. S. is partially supported by ANR project CIMI ANR-11-LABX-0040 and ANR project NQG ANR-23-CE40-0005.

\bibliographystyle{abbrv}
\bibliography{biblio} 

\end{document}